\documentclass{elsarticle}

\biboptions{sort&compress,numbers}
\usepackage[nodots]{numcompress}
\usepackage{lineno,url}
\usepackage{epstopdf}
\usepackage{amsmath,amsfonts,amssymb}
\usepackage{amssymb,amsmath,amsthm}
\usepackage{overpic}
\usepackage{contour}\contourlength{0.7pt}
\usepackage{tikz}\usetikzlibrary{shapes,arrows}
\usepackage{versions}
\usepackage{natbib}
\usepackage{algorithm, algorithmic}
\usepackage{multirow}

\newtheoremstyle{mytheorem}{5pt plus 5pt minus 3pt}{4pt plus 3pt minus 1.5pt}
	{\itshape}{}{\bfseries}{.}{1ex plus 1ex minus .5ex}{}
\newtheoremstyle{mydef}{5pt plus 5pt minus 3pt}{4pt plus 3pt minus 1.5pt}
	{}{0pt}{\bfseries}{.}{1ex plus 1ex minus .5ex}{}
\newtheoremstyle{myremark}{5pt plus 5pt minus 3pt}{4pt plus 3pt minus 1.5pt}
	{}{0pt}{\itshape}{.}{1ex plus 1ex minus .5ex}{}
\theoremstyle{mytheorem}
\newtheorem{prop}{Proposition}[section]

\theoremstyle{mydef}

\newtheorem{ex}[prop]{Example}
\theoremstyle{myremark}
\newtheorem{rem}{Remark}

\newcommand{\de}{\mathrm{d}}

\newcommand{\bm}{\mathbf m}
\newcommand{\bx}{\mathbf x}

\newcommand{\DDt}{\widetilde{\mathbf D}}
\newcommand{\Q}{\mathcal Q}
\newcommand{\Qt}{\widetilde{\Q}}
\newcommand{\taut}{\widetilde{\tau}}
\newcommand{\omegat}{\widetilde{\omega}}
\newcommand{\tf}{\footnotesize}
\newcommand{\br}{\mathbf r}
\def\R {\mathbb R}
\newcommand{\FF}{\mathcal F}
\newcommand{\FFt}{\widetilde{\mathcal F}}

\newcommand{\Om}{\Omega}
\newcommand{\Omt}{\widetilde{\Omega}}

\newcommand{\XX}{\mathcal X}
\newcommand{\XXt}{\widetilde{\mathcal X}}
\newcommand{\xxt}{\widetilde{\bx}}
\newcommand{\bmt}{\widetilde{\bm}}
\newcommand{\mt}{\widetilde{m}}
\newcommand{\xt}{\widetilde{x}}
\newcommand{\St}{\widetilde{S}}
\newcommand{\Dt}{\widetilde{D}}
\newcommand{\Bt}{\widetilde{B}}

\newcommand{\sg}{\hspace{0.1cm}}

\def\Acaption#1#2{\caption{#2}\vspace*{-#1}}

\definecolor{Mgreen}{RGB}{34,139,34}
\definecolor{blau}{rgb}{0.15,0.2,0.5}
\definecolor{gray}{rgb}{0.5,0.5,0.5}
\definecolor{drot}{rgb}{0.7,0,0.1}
\definecolor{gelb}{rgb}{.55,.40,.1}
\definecolor{magenta}{rgb}{1.,0.,1.}
\definecolor{cyan}{rgb}{0.,1.,1.}
\definecolor{green}{rgb}{0.,1.,0.}
\definecolor{Morange}{rgb}{1.,0.5,0.}

\usepackage{array}
\newcolumntype{C}[1]{>{\centering\let\newline\\\arraybackslash\hspace{0pt}}m{#1}}

\begin{document}

\begin{frontmatter}

\title{
Gauss-Galerkin quadrature rules\\ for quadratic and cubic spline spaces\\ and their application to isogeometric analysis}
\author[NumPor]{Michael Barto\v{n}\corref{cor1}}
\ead{Michael.Barton@kaust.edu.sa}
\author[Curtin,AMCS]{Victor Manuel Calo}
\ead{Victor.Calo@Curtin.edu.au}

\cortext[cor1]{Corresponding author}

\address[NumPor]{Center for Numerical Porous Media, King Abdullah University of Science and Technology, Thuwal 23955-6900, KSA}%
\address[AMCS]{Applied Mathematics $\&$ Computational Science and Earth Science $\&$ Engineering, \\
King Abdullah University of Science and Technology, Thuwal 23955-6900, KSA}
\address[Curtin]{Chair in Computational Geoscience Western Australian School of Mines,\\
Faculty of Science and Engineering, Curtin University,\\
Kent Street, Bentley, Perth, Western Australia, 6102, Australia}

\begin{abstract}
We introduce Gaussian quadrature rules for spline spaces that are frequently used
in Galerkin discretizations to build mass and stiffness matrices. By definition,
these spaces are of even degrees. The optimal quadrature rules we recently derived \cite{StiffnessMatrix-2015}
act on spaces of the smallest odd degrees and, therefore, are still slightly sub-optimal.
In this work, we derive optimal rules directly
for even-degree spaces and therefore further improve our recent result.
We use optimal quadrature rules for spaces over two elements as elementary building blocks and
use recursively the homotopy continuation concept described in \cite{Homotopy-2015} to derive optimal rules
for arbitrary admissible number of elements.
We demonstrate the proposed methodology on relevant examples, where we derive optimal
rules for various even-degree spline spaces.
We also discuss convergence of our rules to their asymptotic counterparts, these are
the analogues of the midpoint rule of Hughes et al. \cite{Hughes-2010}, that are exact
and optimal for infinite domains.
\end{abstract}

\begin{keyword}
optimal quadrature rules, Galerkin method, Gaussian quadrature, B-splines, isogeometric analysis, homotopy continuation for quadrature
\end{keyword}

\end{frontmatter}

\section{Introduction and motivation}\label{intro}

Numerical integration is a fundamental ingredient of isogeometric analysis (IGA) and finite elements (FE),
and its computational efficiency is essential.
When simulating physical processes, e.g., \cite{OptimizedWings-2014,DynEarthSol2D-2013,vignal2015coupling,cortes2015performance,lipton2010robustness},
with Galerkin isogeometric discretizations,
specific spline spaces appear when building mass and stiffness matrices. By construction, these spline spaces are of even degrees
and the quadrature rules used to numerically integrate functions from these spaces are sub-optimal \cite{Hughes-2010,Hughes-2012,Hughes-2014,StiffnessMatrix-2015}.
In this work, we present two families of optimal rules for spaces arising from Galerkin discretizations when the original spline degrees are quadratic and cubic.
The derived rules are optimal in terms of number of quadrature points and therefore they minimize the computational cost of the assembly of mass stiffness matrices, 
while guaranteeing exactness of all integrands up to machine precision.
We name these rules \emph{Gauss-Galerkin}.


A \emph{quadrature rule}, or simply a \emph{quadrature}, is an \textit{$m$-point rule},
if $m$ evaluations of a function $f$ are needed to approximate its weighted integral
over a closed interval $[a,b]$
\begin{equation}\label{eq:GaussQuad}
\int_a^b w(x) f(x) \, \mathrm{d}x = \sum_{i=1}^{m} \omega_i f(\tau_i) + R_{m}(f),
\end{equation}
where $w$ is a fixed non-negative \emph{weight function} defined over $[a,b]$.
The rule is required to be \emph{exact}, i.e., $R_m(f) \equiv 0$
for each element of a spline space $S$.
A rule is \emph{optimal} if $m$ is the minimum number of
\emph{weights} $\omega_i$ and \textit{nodes} $\tau_i$ (points at which $f$ is evaluated).

For polynomials, the optimal rule is known to be the classical Gaussian quadrature \cite{Gautschi-1997}
with the \emph{order of exactness} $2m-1$, that is, only $m$ evaluations are needed to exactly integrate
any polynomial of degree at most $2m-1$. Consider a sequence of polynomials $(q_0,q_1,\ldots,q_m,\ldots )$  that form an orthogonal basis
with respect to the scalar product
\begin{equation}
<f,g> = \int_a^b f(x)g(x)w(x) \de x.
\end{equation}
The quadrature points are the roots of the $m$-th orthogonal polynomial $q_m$ which
in the case when $w(x)\equiv 1$ is the degree-$m$ Legendre polynomial \cite{Szego-1936}.

Quadrature rules for polynomial spaces of even degree, however, contain a certain sub-optimality.
For these spaces, even degree implies odd dimension of the space, and therefore $2m-1$ basis functions are being exactly
integrated by $m$ nodes. That is, the rule with $m$ nodes and $m$ weights does not possess
``double precision'', see \cite{Micchelli-1977}, because it integrates (using $2m$ numbers) only a space of dimension $2m-1$.

Looking at this scenario from the point of view of algebraic systems, there are $2m-1$ algebraic equations (constraints that
the quadrature rule exactly integrates the basis functions), but $2m$ unknowns (nodes and weights).
Such an algebraic system is undetermined by one variable and solving it requires a special treatment \cite{SolverUnivar-2011}.
Moreover, the Gaussian quadrature is not unique anymore, in fact, the optimal Gaussian rules form a one-parameter family.
This non-uniqueness, considering a quadrature rule as a zero of a particular algebraic system, makes the rule difficult to
be traced numerically. To make the system well-constrained, one can prescribe one node (weight), e.g., by
considering Gauss-Radau rule where one boundary point is assumed to be a node in the quadrature rule \cite{Gautschi-1997}.

The situation for even-degree spaces is different in the case of splines when compared to polynomials.
The optimal quadrature rules for spline were studied in the late 50's \cite{Schoenberg-1958,Micchelli-1972,Micchelli-1977}.
Micchelli and Pinkus \cite{Micchelli-1977} considered boundary constraints and proved that,
for spaces with uniform continuities (knot multiplicities),
there always exists an optimal quadrature formula with the following number of necessary evaluations:
\begin{equation}\label{eq:Micchelli}
     d + 1 + i = 2 m,
\end{equation}
where $d$ is the polynomial degree, $i$ is the total number of interior knots (when counting multiplicities),
and $m$ is the number of optimal nodes. Moreover, for the case with no boundary constraints,
the quadrature formula is unique, whenever the count in (\ref{eq:Micchelli}) admits it.
Therefore, one can seek optimal rules for even-degree spaces that are unique, contrarily to the polynomial case.

Quadrature rules for splines are important tools for the isogeometric analysis community \cite{Hughes-2010,Hughes-2012,Hughes-2014}
because they are cheap and elegant alternatives to symbolic integration \cite{Gautschi-1997}.
Recently, alternative methods of building mass and stiffness matrices have been proposed \cite{Bert2014,Bert2015,Adel-2015}. They exploit the observation that,
under certain conditions, the optimal convergence rate of the liner system can be achieved despite the fact that the integration rule
is not exact. In this work, however, we focus on quadrature rules that are \emph{exact}, that is, the rules reproduce the integrals 
under affine mappings exactly up to machine precision.

For spline spaces of various degrees and continuities, the optimal (Gaussian) rules were not known until recently \cite{StiffnessMatrix-2015}.
We showed the connection underlying two different optimal rules via continuous transformations between
the corresponding spline spaces.
%
%
Nonetheless, optimal rules for spaces of odd degrees only have been derived because for those spaces the polynomial analogy,
the classical Gaussian quadrature, offers a \emph{unique} and \emph{optimal} counterpart.
Here, we focus on even-degree spline spaces, even though there is no unique optimal polynomial analogy for them, and
show that the right homotopic setup leads to optimal rules for these spaces.
We seek \emph{optimal} quadrature rules for \emph{finite} domains. When these requirements are violated,
we refer the reader to \cite{Hughes-2012} where sub-optimal rules are derived by considering redundantly many nodes, or to
\cite{Hughes-2010} for the case where the integration domain is the whole real line.

Regarding concrete optimal quadrature rules for spline spaces over finite domains, to the best of our knowledge, there are few reports 
in the literature, particularly for even-degree spline spaces. Except for the quadratic case \cite{Hughes-2014}, we are not aware of any
rules for even-degree spline spaces. For odd degrees, we refer the reader to \cite{Nikolov-1996,Quadrature31-2014,Quadrature51-2014}. 
These rules act in turn on uniform cubic, non-uniform cubic, and uniform quintic spline spaces, respectively. These rules are explicit, that is,
there is a formula that computes the node locations and their weights in the first boundary element and a recursion relation derives the other nodes,
parsing from the boundary to the middle of the interval.

For other higher degrees (odd or even), we are not aware of the existence of explicit rules, and therefore a
numerical solver is required. However, to initialize the minimum number of nodes
such that the solver finds a \emph{global} minimum is challenging.
Since the problem is highly non-linear, in general, an arbitrary initial guess usually leads
the solver a local minimum. To avoid such a scenario, we have proposed a methodology
that uses continuity arguments and transfers optimal rules from one space (source) to another (target) by
continuously modifying the knot vector \cite{Homotopy-2015}. In such a setting, the optimal rule is thought of as a root of a particular
polynomial system, and the optimal rule from the previous iteration serves as an initial guess for the updated (modified) spline space.
We showed that the right homotopic setting admits situations where the source and target rules require different
numbers of optimal nodes, and derived optimal rules for several odd-degree spaces of various continuities \cite{StiffnessMatrix-2015}.

In this paper, we extend our recent results on optimal quadrature rules
to derive optimal rules for spaces of even degrees.
We show that, using two-element optimal quadrature blocks, one can build an appropriate source space and
set up the homotopy continuation to derive optimal rules for appropriate target spaces.

The rest of the paper is organized as follows. Section~\ref{sec:matrices} describes the particular class of spline spaces
for which the optimal rules are derived. Section~\ref{sec:splines} summarizes a few basic properties of spline spaces
and shows possible continuous transformations between them. In Section~\ref{sec:homotopy}, we recall how
the homotopy continuation is applied to derive new optimal quadrature rules.
Section~\ref{sec:ex} shows the results of the derived optimal rules, discusses their validity, and points the connection
to their asymptotic counterparts.
The paper is concluded by discussing our observations and future research directions in Section~\ref{sec:conl}.

\section{Spline spaces appearing in Galerkin discretizations}\label{sec:matrices}


When solving elliptic partial differential equations using weighted-residual methods such as
finite elements and isogeometric analysis,
one builds Grammian (mass) and stiffness matrices by computing the $L_2$ scalar products
of the $l$-th derivatives of the basis functions. Considering the original spline space of
degree $p$ and continuity $k$, $S_{p,k}$, the mass matrix contains elements from $S_{2p,k}$.
The stiffness matrix, depending on the order $l$ of the differential operator, contains elements
from $S_{2(p-l),k-l}$. The hierarchy of spline spaces is shown in Fig.~\ref{fig:hierarchy}.
Even though one could use different quadrature rules for each of these spaces, common
practice is to minimize the computational cost and simplify implementation
by using one rule that integrates exactly elements from both spaces \cite{petiga}. The minimum
spline space that contains both of them is $S_{2p,k-l}$ and therefore we focus on this category of spaces
of degree $d:=2p$ and continuity $c:=k-l$ for various $p$, $k$, and $l$, and derive optimal quadrature rules for them.
To the best of our knowledge, there are no optimal rules derived for $S_{2p,k-l}$, $p>1$.

\begin{rem}
In 2D and 3D, other intermediate products need to be integrated, but these still belong to $S_{2p,k-l}$.
Thus, these are not discussed further.
\end{rem}

 \begin{figure}[!tb]
\vrule width0pt\hfill
\begin{equation*}
\renewcommand{\arraystretch}{1.4}
\begin{array}{ccc}
(p,k)     & \rightarrow & (2p,k)   \\
(p-1,k-1) & \rightarrow & (2(p-1),k-1)   \\
\vdots    & & \vdots \\
(p-l,k-l) & \rightarrow & (2(p-l),k-l)
\end{array}
\end{equation*}
 \hfill \vrule width0pt\\
 \vspace{-25pt}
\Acaption{1ex}{A hierarchy of spline spaces used when building mass and stiffness matrices.
Left: the original spline space of degree $p$ and continuity $k$ is being differentiated $l$-times.
Right: the corresponding spline spaces that contain scalar products are shown.}\label{fig:hierarchy}
 \end{figure}


\begin{ex}
Let $p=3$, $k=2$, $l=1$, then we have
\begin{equation*}
\renewcommand{\arraystretch}{1.4}
\begin{array}{ccc}
(6,2)     & \subset & (6,1)   \\
\cup      &  & \cup   \\
(4,2)     & \subset & (4,1)
\end{array}
\end{equation*}
The inclusion relations follow directly from the fact that the corresponding knot vectors
are nested \cite{deBoor-1978}. In the context of finite elements in 1D,
the elements of the mass matrix belong to a $(6,2)$-space
while the scalar products that fill the stiffness matrix belong to a $(4,1)$-space.
Thus, we seek the smallest spline space that contains both
$(6,2)$ and $(4,1)$, i.e, $(d,c)=(6,1)$.
\end{ex}



\section{Continuous transformations between spline spaces}\label{sec:splines}

Consider a spline space defined above a knot vector
\begin{equation}\label{eq:x}
\begin{array}{ccccccc}
\XX_N =  & (a= &\underbrace{x_0, \dots,x_0,} & \underbrace{x_1, \dots, x_{1},} & \dots & \underbrace{x_N,\dots,x_N}& =b)\\
& & m_0 & m_1 & & m_N &
\end{array}
\end{equation}
which, for the sake of simplicity, we split $\XX_N:=(\bx,\bm)$ into the \emph{domain partition} $\bx$, $\bx \in \mathbb R^N$,
and the \emph{vector of multiplicities} $\bm$, $\bm \in \mathbb N^N$, and write
\begin{equation}\label{eq:split}
\bx =  (x_0, \dots,x_N), \quad \quad \bm = (m_0,\dots, m_N).
\end{equation}
We further recall $1\leq m_i \leq d+1$, $i=0,\dots,N$ and assume $\XX_N$ is an \emph{open knot} vector on $[a,b]$, that is, $m_0=m_N = d+1$,
and $d$ is even.
We denote by $\pi_d$ a space of polynomials of degree at most $d$ and define
the spline space associated to $\XX_N$ as
\begin{equation}\label{eq:targetSpace}
S_{\bx,\bm}^{N,d} = \{ f\in C^{d-m_k} \,\, \textnormal{at} \, x_k, k=0,\dots,N \,\, \textnormal{and} \,  f|_{(x_{k-1},x_{k})} \in \pi_d, k=1,\dots,N\}.
\end{equation}
Our goal is to derive a Gaussian rule for this target space $S_{\bx,\bm}^{N,d}$.
To do so, we define an associated source space for which the optimal rule is known. 
Consider a \emph{source} knot vector
\begin{equation}\label{eq:xt}
\begin{array}{ccccccc}
\XXt_n =  & (a= &\underbrace{\xt_0, \dots,\xt_0,} & \underbrace{\xt_1, \dots, \xt_{1}},& \dots & \underbrace{\xt_n,\dots,\xt_n}& =b)\\
& & \mt_0 & \mt_1 & & \mt_n &
\end{array}
\end{equation}
and using analogous notation to (\ref{eq:split}), i.e., $\XXt_n:=(\xxt,\bmt)$, we obtain
\begin{equation}\label{eq:splitt}
\xxt =  (\xt_0, \dots, \xt_n), \quad \bmt = (\mt_0,\dots, \mt_n),
\end{equation}
where $\xt$ is an arbitrary (uniform or non-uniform) partition and $\bmt$ are the corresponding multiplicities.
%
We define the source spline space as
\begin{equation}\label{eq:sourceSpace}
\St_{\xxt,\bmt}^{n,d} = \{ f\in C^{-1} \,\, \textnormal{at} \, \xt_k, k=0,\dots,n \,\, \textnormal{and} \,  f|_{(\xt_{k-1},\xt_{k})} \in \pi_d, k=1,\dots,n\}.
\end{equation}
Our aim is to continuously transform $\St_{\xxt,\bmt}^{n,d}$ into $S_{\bx,\bm}^{N,d}$ over $[a,b]$,
controlling the transformation by a continuous knot transformation
\begin{equation}\label{eq:Trans}
\XXt_n \rightarrow \XX_N.
\end{equation}
%
%

\begin{figure}[!tb]
\vrule width0pt\hfill
    \begin{overpic}[width=.89\textwidth,angle=0]{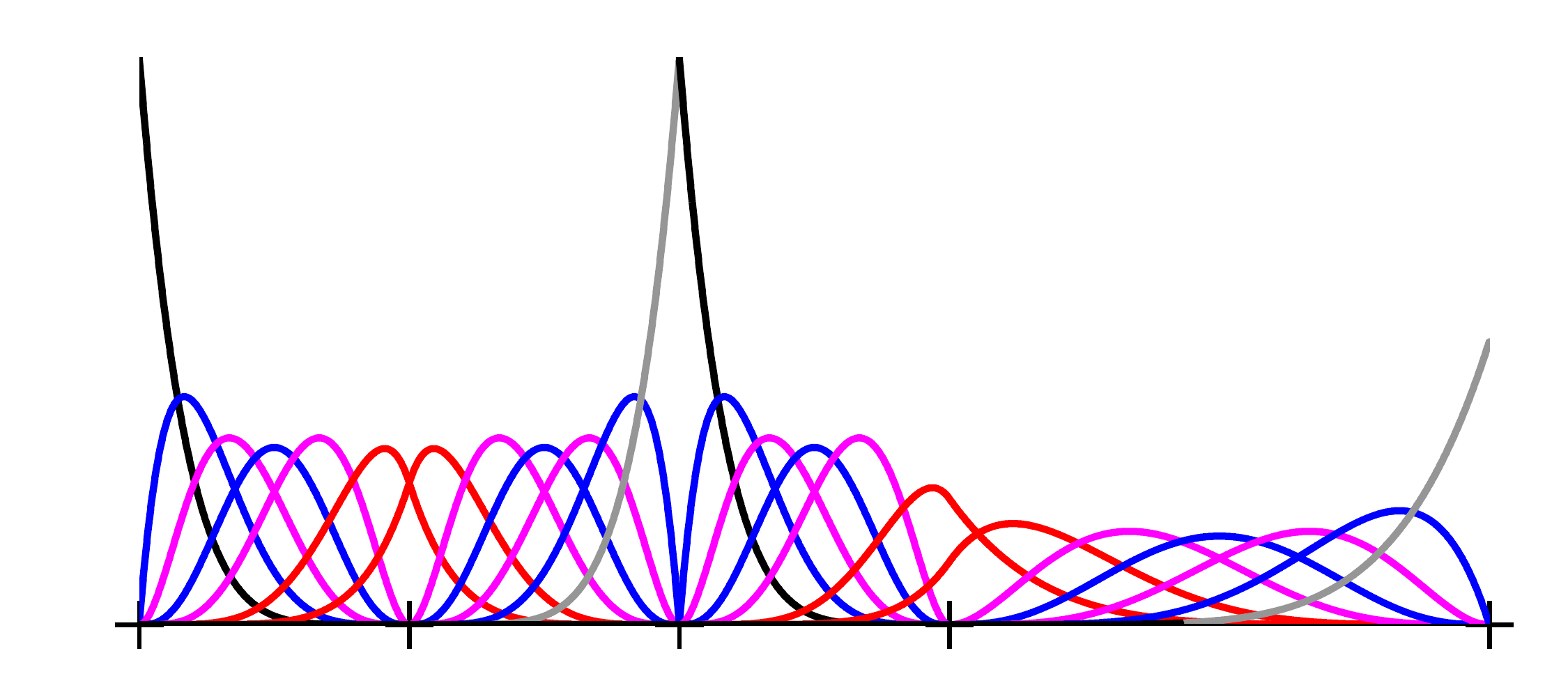}
        \put(11,35){\small$\Dt_{2d(k-3)+1}$}
        \put(25,0){\small $\xt_{k-2}$}
        \put(4,0){\small$\xt_{k-3}$}
        \put(95,0){\small$\xt_{k+1}$}
        \put(43,0){\small$\xt_{k-1}$}
        \put(58,0){\small$\xt_{k}$}
        \put(85,25){\small$\Dt_{2d(k-1)+12}$}
    \end{overpic}
\hfill \vrule width0pt\\[-2ex]
\Acaption{1ex}{A part of the source spline space of degree six over four elements is generated from two basic blocks, each of them consisting of two elements.
The blocks meet with $C^{-1}$-continuity at $\xt_{k-1}$, while the continuity inside the blocks is $C^1$ ($\xt_{k-2}$ and $\xt_k$ have multiplicity five).
 }\label{fig:SourceSpace61}
 \end{figure}

As an example, consider a source space that consists of $n/2$ two-element blocks with $C^1$-continuity inside the block and $C^{-1}$-continuity
at the blocks' boundary. For this particular type of source space $\St^{n,d}_{\xxt,\bmt}$,
we generate the basis $\DDt = \{\Dt_i\}_{i=1}^{nd}$ as
\begin{equation}\label{eq:Bernstein}
\renewcommand{\arraystretch}{1.2}
\begin{array}{lcl}
\Dt_{2d(k-1)+1}(t)   & = & [\xt_{k-1},\xt_{k-1},\dots, \xt_{k-1},\xt_k](. - t)_{+}^{d}   \\
 &\vdots & \\
\Dt_{2d(k-1)+d}(t)   & = & [\xt_{k-1},\xt_{k-1},\xt_{k},\dots,\xt_k,\xt_{k+1}](. - t)_{+}^{d}   \\
\Dt_{2d(k-1)+d+1}(t) & = & [\xt_{k-1},\xt_{k},\dots,\xt_k,\xt_{k+1},\xt_{k+1}](. - t)_{+}^{d}   \\
 &\vdots & \\
\Dt_{2d(k-1)+2d}(t)   & = & [\xt_{k},\xt_{k+1},\dots, \xt_{k+1},\xt_{k+1}](. - t)_{+}^{d}, \\
\end{array}
\end{equation}
where $[.]f$ stands for the divided difference, $u_{+} = \max(u,0)$ is the truncated
power function, and $k=1,\dots, n/2$. An example of a source space is shown in Fig.~\ref{fig:SourceSpace61}.
We work with non-normalized basis and therefore
\begin{equation}\label{interiorIntegral}
I[\Dt_i] = \frac{1}{d+1}\; \textnormal{for} \quad i = 1,2,\ldots,nd,
\end{equation}
where $I[f]$ stands for the integral of $f$ over the interval $[a,b]$, see e.g. \cite{Hoschek-2002-CAGD}.

We build our source space $\St_{\xxt,\bmt}^{n,d}$ for which an optimal quadrature rule is known.
In particular, we use a union of spline spaces, \emph{building blocks}. The blocks meet with $C^{-1}$ continuity, and each
block is provided with an optimal quadrature rule. 
Since these blocks have optimal quadrature rules and, due to the $C^{-1}$ continuity between these blocks, the union
space also has an optimal rule.
We generate these spaces recursively, starting with a single block containing \emph{two} elements which are connected with
the desired continuity. 
For example, for $d=6$, $c=1$, $n=2$, the optimal two-element rule requires six nodes as \eqref{eq:Micchelli} becomes $7+1+5=12$.
We now detail the construction of the initial building blocks.

\subsection{Two-element building blocks}\label{ssec:blocks}

The difficulty of optimal rules for even degree spaces stems from the fact that the classical polynomial
Gauss quadrature for even degrees is suboptimal in the sense that $m$ quadrature points are used
for $2m-1$ basis functions. In fact, there is a whole one parameter family of Gaussian rules.
Using such a rule as our source rule, the sub-optimality grows with the number of elements and, therefore, starting the homotopy
continuation algorithm with a highly sub-optimal initial rule would require vanishing too many nodes from the system.

\begin{rem}
An alternative to a one parameter family of Gauss rules for even degree polynomials is to select a node location
a priori to make the polynomial system well-constrained. In particular, Gauss-Radau rules \cite{Gautschi-1997} choose
one of the endpoints as a node.
\end{rem}

%

\begin{figure}[!tb]
\vrule width0pt\hfill
    \begin{overpic}[width=.5\textwidth,angle=0]{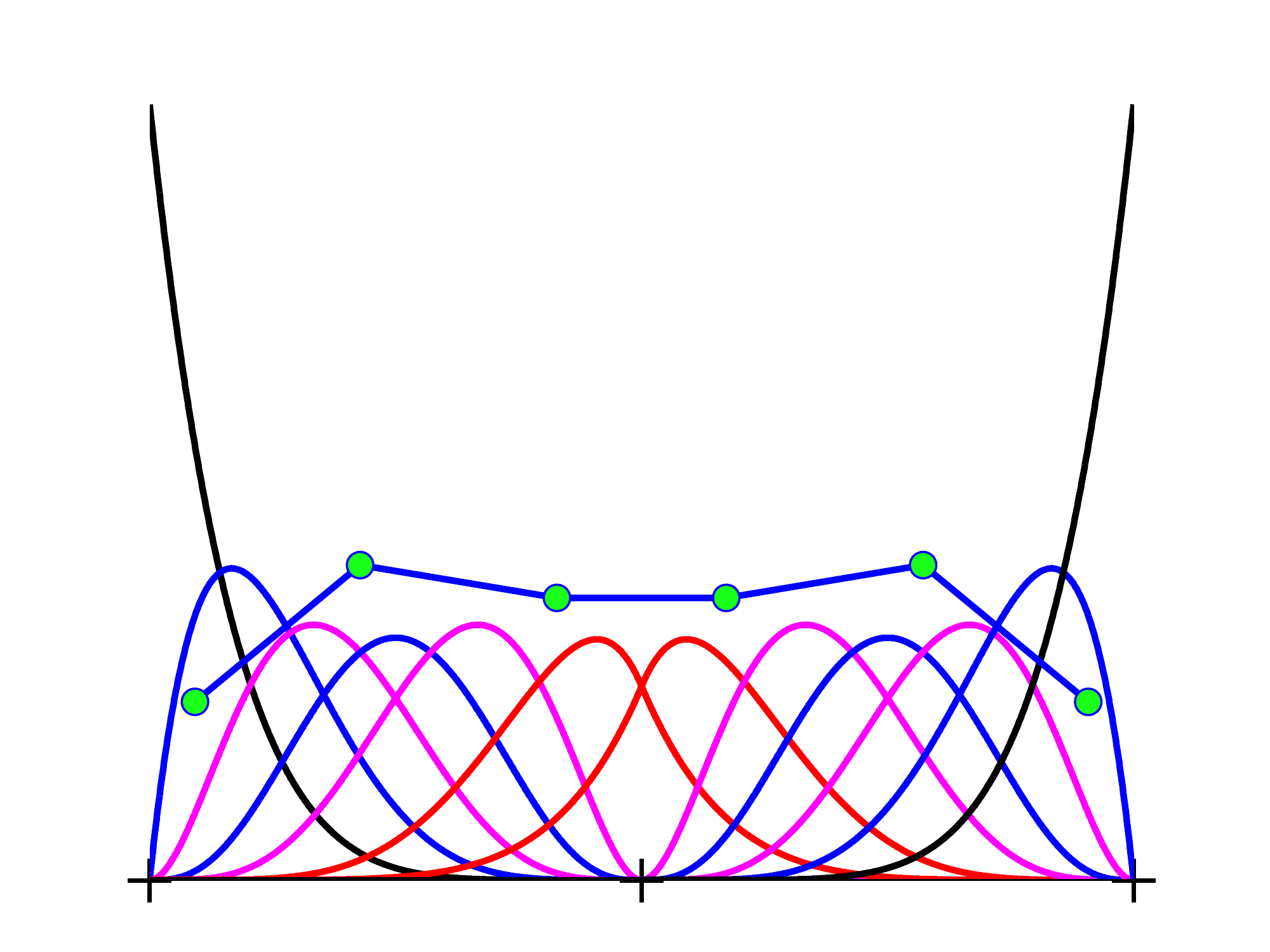}
        \put(15,55){\small$\Dt_1$}
        \put(100,60){\small$\dim(\St^{n,d}_{\xxt,\bmt})=12$}
        \put(100,50){\small $\xxt = (x_0,x_1,x_2)$}
        \put(100,40){\small $\bmt = (7,5,7)$}
        \put(25,35){\small $[\tau_2,\omega_2]$}
        \put(5,0){\small$a=\xt_{0}$}
        \put(85,0){\small$b=\xt_{2}$}
        \put(45,0){\small$\xt_{1}$}
    \end{overpic}
\hfill \vrule width0pt\\[-2ex]
\Acaption{1ex}{Gauss type initial building block. Optimal quadrature rule over a two-element ($n=2$) block for $d=6$, $c=1$.
The above-built spline space $\St^{n,d}_{\xxt,\bmt}$ is of dimension $12$.
Due to the $C^1$ continuity at $\xt_{1}$, two basis functions span both elements (red). The optimal
quadrature rule requires six quadrature points (green dots) and is computed, due to symmetry, from the system \eqref{eq:IniSystem61a}.
 }\label{fig:BasicBlock61}
 \end{figure}

Depending on the dimension of the space, 
we build a two-element source space provided with a Gauss rule (even dimension) or Gauss-Radau (odd).
As an example, consider $d=6$, $c=1$, $n=2$, see Fig.~\ref{fig:BasicBlock61}. To derive the optimal rule for this particular
spline space, we build the corresponding algebraic system which represents the constraint that the quadrature must integrate exactly
the basis of the space. The dimension of the spline space is twelve, however, due to symmetry,
the polynomial system to solve is $6\times6$
\begin{equation}\label{eq:IniSystem61}
\Qt_a^b[\Dt_i] = I[\Dt_i], \quad i = 1,\dots, 6
\end{equation}
where the source quadrature rule $\Qt$ (linear operator) is expressed in terms of six unknowns $\tau_j, \omega_j$ for $j=1,2,3$.
The particular setting on $[a,b] = [0,2]$ gives
\begin{equation}\label{eq:IniSystem61a}
\renewcommand{\arraystretch}{1.4}
\begin{array}{rcrcrcc}
             \omega_1 (1-\tau_1)^6 & + &            \omega_2 (1-\tau_2)^6 & + &            \omega_3 (1-\tau_3)^6 & = & \frac{1}{7},  \\
  6 \tau_1   \omega_1 (1-\tau_1)^5 & + & 6 \tau_2   \omega_2 (1-\tau_2)^5 & + & 6 \tau_3   \omega_3 (1-\tau_3)^5 & = & \frac{1}{7},  \\
 15 \tau_1^2 \omega_1 (1-\tau_1)^4 & + &15 \tau_2^2 \omega_2 (1-\tau_2)^4 & + &15 \tau_3^2 \omega_3 (1-\tau_3)^4 & = & \frac{1}{7},  \\
 20 \tau_1^3 \omega_1 (1-\tau_1)^3 & + &20 \tau_2^3 \omega_2 (1-\tau_2)^3 & + &20 \tau_3^3 \omega_3 (1-\tau_3)^3 & = & \frac{1}{7},  \\
 15 \tau_1^4 \omega_1 (1-\tau_1)^2 & + &15 \tau_2^4 \omega_2 (1-\tau_2)^2 & + &15 \tau_3^4 \omega_3 (1-\tau_3)^2 & = & \frac{1}{7},  \\
  6 \omega_1 \tau_1^5 -5\omega_1 \tau_1^6 & + & 6 \omega_2 \tau_2^5 -5\omega_2 \tau_2^6 & + & 6 \omega_3 \tau_3^5 -5\omega_3 \tau_3^6 & = & \frac{1}{7},
\end{array}
\end{equation}
and using computer algebra, the sequential factorization yields a univariate polynomial
\begin{equation}\label{eq:Uni61}
2-54\tau_1 + 507\tau_1^2-2024\tau_1^3+3840\tau_1^4-3402\tau_1^5+1127\tau_1^6
\end{equation}
which is solved numerically. There are five real roots in $[0,1]$, two of them being false positive answers.
The reason for that is that symbolic factorization is based on resultant computation, 
and this may introduce artificial roots (false positives of the initial system).
Finally we obtain
\begin{equation}\label{eq:NW61}
\renewcommand{\arraystretch}{1.1}
\begin{array}{ccccccc}
\tau_1    & = & 0.21132486540518711775, & & \omega_1 & = & 0.23004836288935413032 \\
\tau_2    & = & 0.42759570120004222829, & & \omega_2 & = & 0.40614522687566702979  \\
\tau_3    & = & 0.82792440129801198117, & & \omega_3 & = & 0.36380641023497883991
\end{array}
\end{equation}
which determines the optimal quadrature rule for the two-element block.

Depending on the parity of dimension of the spline space, the corresponding type of the optimal initial block rule (Gauss or Gauss-Radau)
is chosen, see Table~\ref{tab:Initial}.


\begin{table}[!tb]
 \begin{center}
  \begin{minipage}{0.9\textwidth}
\caption{Initial two-element block rules for Galerkin spaces of degree $d$ and continuity $c$.
Dimension of the open end knot spline spaces, number of nodes, and the type of the rule (Gauss or Gauss-Radau) are shown.
}\label{tab:Initial}
  \end{minipage}
\vspace{0.2cm}\\
\small{
\renewcommand{\arraystretch}{1.15}
\renewcommand{\tf}{\small}
\begin{tabular}{| r | c | c| c|}\hline
$(d,c)$ & dim & $\#$nodes & type \\\hline\hline
(4,0)   & 9  & 5 & G-R   \\\hline
(6,1)   & 12 & 6 & G   \\\hline
\end{tabular}
}
\end{center}
\end{table}


\section{Gaussian quadrature via homotopy continuation}\label{sec:homotopy}

In this section, we derive optimal quadrature rules for spline spaces with various polynomial degrees and continuities.
We use homotopy continuation as recently introduced in \cite{Homotopy-2015}
and refer the reader to that work for a more detailed description of the methodology.
For the sake of convenience, however, we recall the main ideas used in the framework.


\subsection{Gaussian quadrature}\label{ssec:gauss}

We consider a source space $\St_{\xxt,\bmt}^{n,d}$ over $n$ elements
 with an optimal quadrature formula
\begin{equation}\label{quadratureS}
\Qt_a^b[f] = \sum_{i=1}^{m} \omegat_i f(\taut_i) = \int_{a}^{b} f(t) \mathrm{d}t, \quad f \in \St_{\xxt,\bmt}^{n,d}.
\end{equation}
and consider $\St_{\xxt,\bmt}^{n,d}$ as a union of two (or more)
spline spaces (elementary blocks) connected with $C^{-1}$-continuity such that each block is provided
with an optimal rule.

Consider the target space $S_{\bx,\bm}^{N,d}$ and let $r$
be an even dimension difference
between the source and the target spaces. Then, according to \cite{Micchelli-1977},
there exists an optimal target rule
\begin{equation}\label{quadrature}
\Q_a^b[f] = \sum_{i=1}^{m-\frac{r}{2}} \omega_i f(\tau_i)= \int_{a}^{b} f(t) \mathrm{d}t, \quad f \in S_{\bx,\bm}^{N,d}.
\end{equation}

The source and the target rules do not require the same number of optimal quadrature points and the transition
of the rule is still possible, we refer the reader to \cite[Section 4]{StiffnessMatrix-2015}.
The proposed methodology, see \cite{Homotopy-2015}, traces the quadrature rule,
as the source space is transformed into the target space by continuously modifying the knot vector, via (\ref{eq:Trans}).
The quadrature rule, $\Q$, represented by its nodes and weights, is a function of time $t$, $t\in [0,1]$.
If no ambiguity is imminent, we omit the time parameter and write $\tau_i$ instead of $\tau_i(t)$.
The source rule is $\Qt=\Q(0)$ and the target rule we wish to derive is $\Q=\Q(1)$.

\begin{figure}[!tb]
\vrule width0pt\hfill
 \begin{overpic}[width=0.89\textwidth,angle=0]{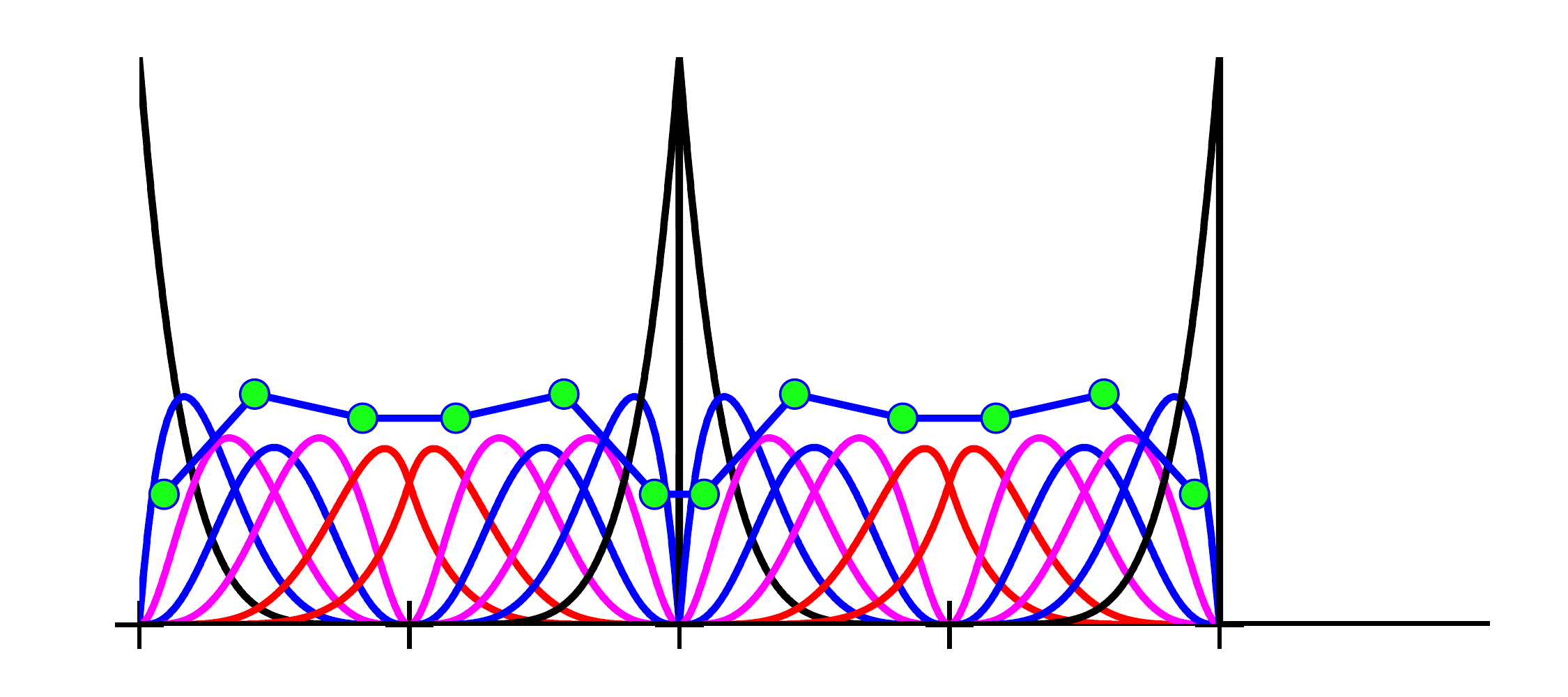}
    \put(12,35){\fcolorbox{gray}{white}{$\St_{\xxt,\bmt}^{4,6}$}}
    \put(50,28){\scriptsize$\xxt = (0,1,2,3,4)$}
    \put(50,32){\scriptsize$\bmt = (7,5,7,5,7)$}
    \put(50,36){\scriptsize$\dim(\St_{\xxt,\bmt}^{4,6}) = 24$}
    \put(7,6){\scriptsize$a$}
    \put(79,6){\scriptsize$b$}
    \put(9,1){\scriptsize$7$}
    \put(25,1){\scriptsize$5$}
    \put(42,1){\scriptsize$7$}
    \put(61,1){\scriptsize$5$}
    \put(76,1){\scriptsize$7$}
    \put(25,22){$\Bt_0^L$}
    \put(60,22){$\Bt_0^R$}
    \end{overpic}\hfill \vrule width0pt\\[-1ex]
\vrule width0pt\hfill
 \begin{overpic}[width=.89\columnwidth,angle=0]{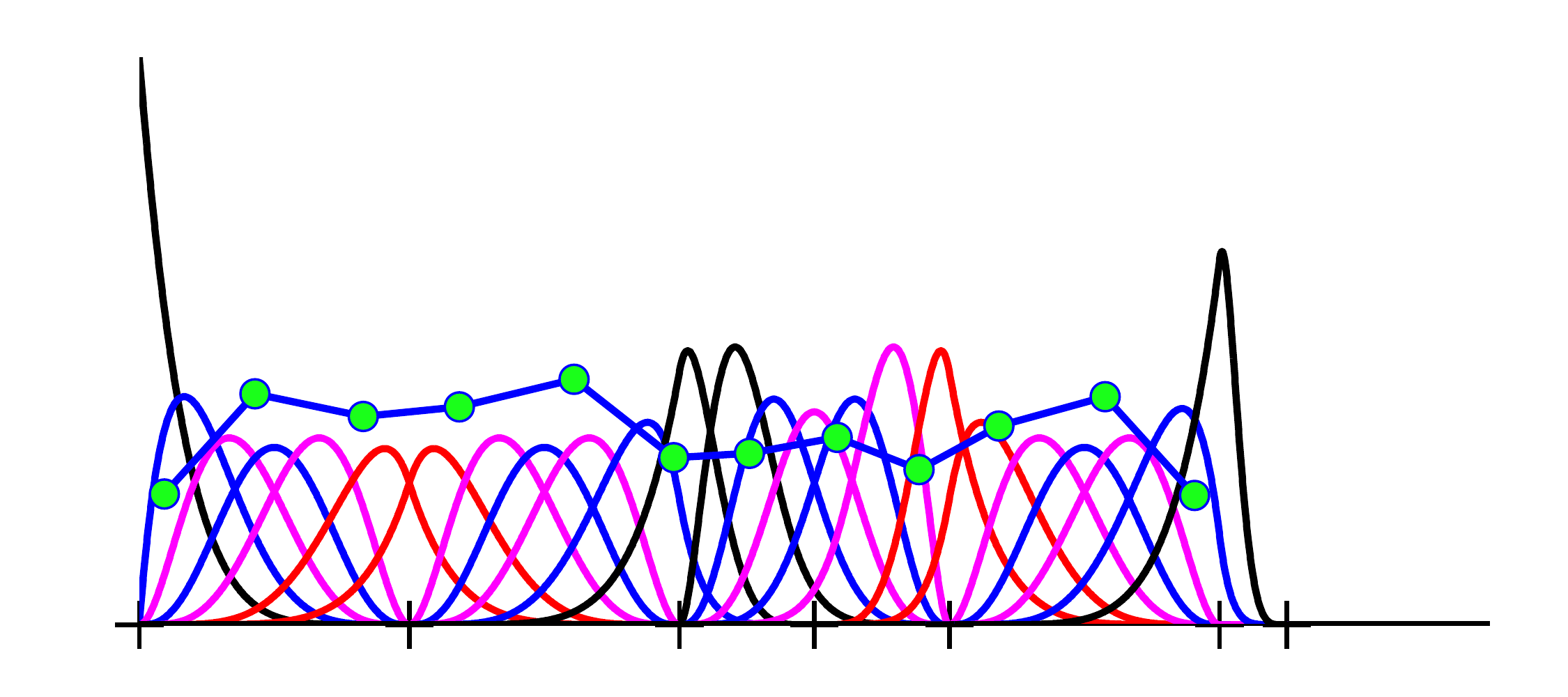}
    \put(9,1){\scriptsize$7$}
    \put(25,1){\scriptsize$5$}
    \put(42,1){\scriptsize$5$}
    \put(61,1){\scriptsize$5$}
    \put(76,1){\scriptsize$5$}
    \put(52,1){\scriptsize$2$}
    \put(82,1){\scriptsize$2$}
	\end{overpic} \hfill \vrule width0pt\\[-1ex]
\vrule width0pt\hfill
 \begin{overpic}[width=.89\columnwidth,angle=0]{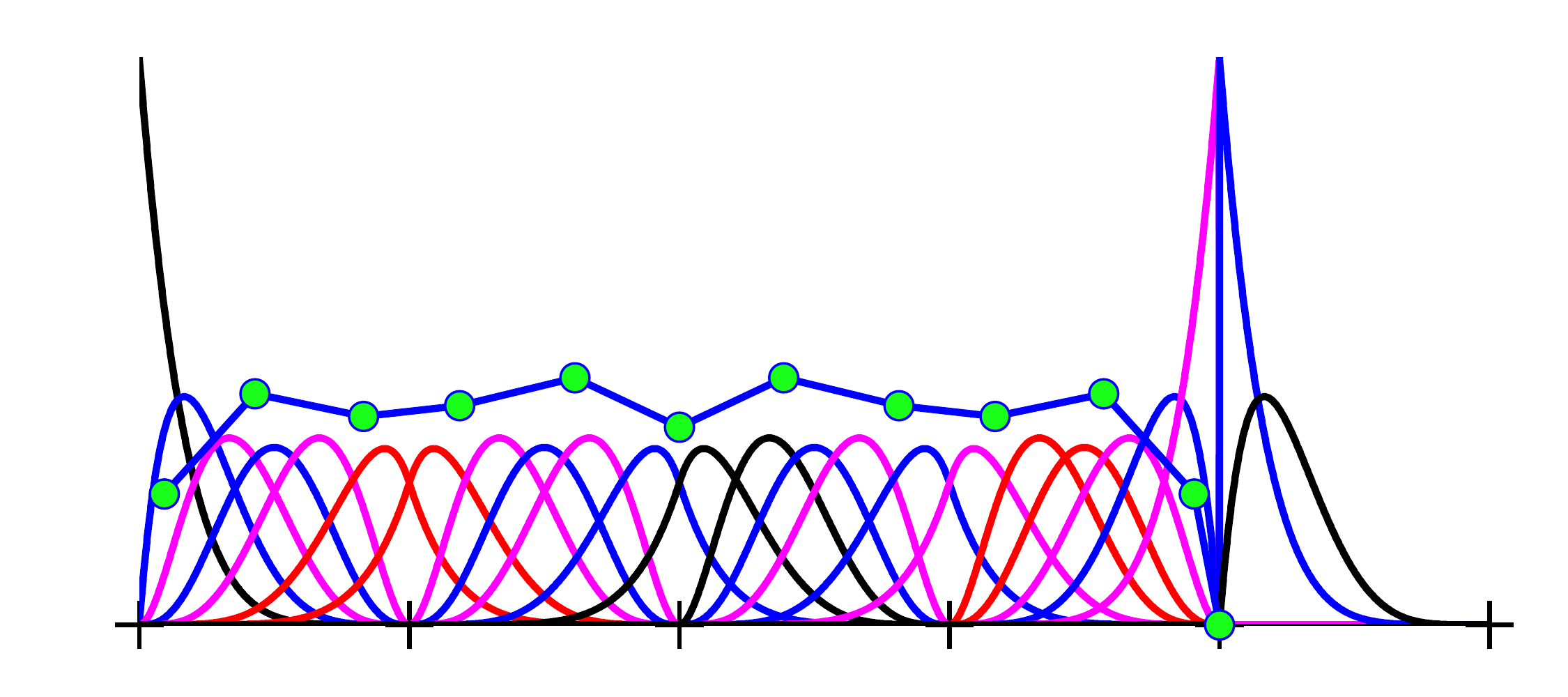}
    \put(12,35){\fcolorbox{gray}{white}{$S_{\bx,\bm}^{4,6}$}}
    \put(50,28){\scriptsize$\bx = (0,1,2,3,4,5)$}
    \put(50,32){\scriptsize$\bm = (7,5,5,5,7,2)$}
    \put(50,36){\scriptsize$\dim(S_{\bx,\bm}^{2,7}) = 14$}
    \put(80,35){\scriptsize$D_{23}$}
    \put(84,16){\scriptsize$D_{24}$}
    \put(79,3){\fcolorbox{gray}{white}{\scriptsize$b=\tau_{12}$}}
    \put(7,6){\scriptsize$a$}
    \put(9,1){\scriptsize$7$}
    \put(25,1){\scriptsize$5$}
    \put(42,1){\scriptsize$5$}
    \put(61,1){\scriptsize$5$}
    \put(76,1){\scriptsize$7$}
    \put(95,1){\scriptsize$2$}
	\end{overpic}\hfill \vrule width0pt\\[-4ex]
\Acaption{1ex}{The continuous evolution of optimal quadrature rules. Top: the initial source space of degree $d=6$, $\St^{4,6}_{\xxt,\bmt}$,
 spans four uniform elements ($n=4$) and is generated from two elementary building blocks $\Bt_0^L$ and $\Bt_0^R$
 for which the optimal source rules are known (green dots). The blocks are merged with $C^{-1}$ continuity;
 vertical lines highlight the discontinuities.
  The integers bellow knots represent their multiplicities. Middle: an intermediate step of the continuation. Two double knots are being moved; 
  one moves from the middle of $[a,b]$ towards $b$ while the second one leaves $[a,b]$. Bottom: the target space $S_{\bx,\bm}^{4,6}$ together with its optimal rule.
  The dimension of $S_{\bx,\bm}^{4,6}$ on $[a,b]$ is $22$ since two basis functions, $D_{23}$ and $D_{24}$, lost their support on $[a,b]$. Consequently
  $\tau_{12}=b$ and $\omega_{12}=0$ which is in accord with the fact that the optimal rule for the target space requires only $11$ nodes, see \eqref{eq:Micchelli}.}\label{fig:Trans61}
 \end{figure}

\subsection{Homotopy continuation}\label{ssec:homotopy}

Polynomial homotopy continuation is a numerical scheme commonly used to solve
polynomial systems of equations \cite{Wampler-1990,Wampler-2005}.
Given a polynomial system $\FF(\bx)=\mathbf 0$ that we want to solve,
the method uses the known roots of a simpler polynomial system (source) $\FFt(\bx)=\mathbf 0$
which is continuously transformed into the desired (target) solution.
We can therefore write
\begin{equation}\label{eq:SystemTime}
\FF(\bx,t)=\mathbf 0
\end{equation}
that at $t=0$ is the system whose roots we know, and at $t=1$ is the target system we aim to solve.
We recommend the reader the book \cite{Wampler-2005}
for a detailed explanation of polynomial homotopy continuation.

In the context of optimal quadrature rules for splines, the homotopic setting is adapted as follows:
the quadrature rule $\Q$ (the nodes and weights) is considered as a point in high-dimensional space
$$
\bx = (\tau_1,\dots, \tau_m, \omega_1,\dots, \omega_m), \quad \bx \in \R^{2m}.
$$
Our source $2m \times 2m$ polynomial system $\FF(\bx,0) = \mathbf 0$ expresses that the source rule $\Qt$,
as stated in equation (\ref{quadratureS}),
exactly integrates the source basis $\DDt$, that is,
\begin{equation}\label{eq:IniSystem}
\Qt_a^b[\Dt_i] = I[\Dt_i], \quad i = 1,\dots, 2m
\end{equation}
and the source root $\br$ that solves (\ref{eq:IniSystem}) is the union of the optimal rules acting on each particular building block
(a subspace of $\St$).

At every instant, a certain domain $\Om \in \R^{2m}$ bounds the root.
For the source domain $\Omt \subset \R^{2m}$ we know the element (knot span) of every node.
For example, for the rule shown in Fig.~\ref{fig:BasicBlock61} we have
\begin{equation}\label{eq:OmegaNodes}
\begin{array}{ccc}
(\tau_1, \dots, \tau_6) & \in & [\xt_0,\xt_1]\times [\xt_0,\xt_1]\times [\xt_0,\xt_1]\times
[\xt_1,\xt_2]\times [\xt_1,\xt_2]\times [\xt_1,\xt_2].
\end{array}
\end{equation}

For the weights we use (a rough) range $[0,b-a]$. Combined together, the source domain is
\begin{equation}\label{eq:OmegaT}
\begin{array}{cccc}
\Omt = & \underbrace{[\xt_0,\xt_1] \times \dots \times [\xt_{n-1}, \xt_n]} & \times & \underbrace{[0,b-a] \times \dots \times [0,b-a]}.\\
& m & & m
\end{array}
\end{equation}
As the source space continuously evolves to the target one, the system $\FF(\bx,t)=\mathbf 0$ continuously changes too,
and so does $\Om(t)$. The root $\br(0)$ of $\FF(\bx,0)= \mathbf 0$ is numerically traced and the root $\br(1)$
of $\FF(\bx,1)= \mathbf 0$ is returned. We refer the reader to \cite{Homotopy-2015} for a detailed description
of this numerical tracing.

In our setup, the knots move towards the right boundary $b$.
In the case of Gauss type rules (see Fig.~\ref{fig:Trans61}), the limit algebraic system (as $\tau_m \rightarrow b$)
is set accordingly by ignoring the last two integral constraints ($I[D_{i}] = \Q_a^b[D_{i}], i=2m-1,2m$). 
The reason for $\tau_m \rightarrow b$ is the fact that two knots (one double knot)
move from inside towards $b$ and the rule must exactly integrate all the basis functions, including the last two $D_{2m-1}$ and $D_{2m}$. 
For example, $D_{2m}$ has a non-zero support on $[x_j,b]$, and since $x_j$ is being moved to $b$, 
$D_{2m}$ looses its support over $[a,b]$.
Consequently, the last two equations of the system become ill-posed ($0=0$) in the limit, and must be removed from the system.
In our implementation, the numerical threshold was set $\varepsilon = I[D_{2m}] = 10^{-3}$.

In the case of Gauss-Radau rules we proceed as follows. To keep the argument simple, let us assume the target space is symmetric with 
respect to the middle of the interval. The dimension of the system is odd and therefore we have one degree of freedom to choose either a node or a weight. 
We require one node to reach a predefined position, the middle of the interval. 
The limit algebraic system is built accordingly by setting $\tau_{\frac{m+1}{2}} = \frac{a+b}{2}$,
remaining only its weight $\omega_{\frac{m+1}{2}}$ as an unknown,
see also later Remark~\ref{rem:GR} and Fig.~\ref{fig:40}. 

\subsection{Building the source space}\label{ssec:source}

We build the source spaces recursively by merging two spline spaces in $C^{-1}$-fashion,
each of these spaces has an optimal rule, see Fig.~\ref{fig:Trans61}.
We use the source space with a known \emph{optimal} rule to apply homotopy
continuation can be applied to transform the source space and to derive the optimal rule for the target space.
In the first iteration, the spline spaces being merged are the elementary building blocks, see also Section~\ref{ssec:blocks}.
In the next iteration, the spline spaces with the computed optimal rules are used to build the source space.

To simplify the arguments, we concentrate on target spaces with uniform continuities, that is, the open end knot vector of multiplicities is
\begin{equation}\label{eq:source}
\bm = (d+1,d-c,\dots, d-c,d+1)
\end{equation}
where $c$ is the desired continuity between elements.

Let $N$, $d$, and $c$ be given. In the zero-th iteration, we use building blocks consisting of two elements
as explained in Section~\ref{ssec:blocks}, and build the source space from them.
We derive the target optimal rule by applying the homotopy continuation, see Fig.~\ref{fig:Trans61},
and progress recursively. 
Let us denote by $\Bt_i^L$ and $\Bt_i^R$ the left and right block spline spaces
that are being merged in the $i$-th level of recursion. Such a merging operation is a union of knot and multiplicities vectors.
Under the assumption that $\Bt_i^L$ and $\Bt_i^R$
have both a known optimal rule, due to the discontinuity,
the optimal rule of the merged space is just a union of the left and right rules.
The sub-spaces to form the source space, $\Bt_i^L$ and $\Bt_i^R$,  need not to have the same dimension.
For example for $N=20$, $d=4$, $c=0$ (Gauss-Radau type rule),
one builds the source space by merging the spaces with the optimal 
rules over four (see Fig.~\ref{fig:SourceSpace40}) and sixteen (see Fig.~\ref{fig:40} top) elements.

\begin{rem}\label{rem:Setup}
Our approach described in this section is a one particular choice of the homotopic setup. One can use different source rules
and different knot transformations to derive the same target rule, see \cite[Section 5]{Homotopy-2015}. We recall the metaphor
of the homotopic evolution of a quadrature rule as a curve in $\mathbb R^{2m}$ connecting two points (the source and target rules). 
Using this metaphor, one can reach the target point using various paths starting from one, or several different, source points. 
Since our objective was to derive the actual rules with a very high precision which we achieved, 
we did not experiment with different source points nor paths.
The knot transformation between two neighboring knot vectors was set uniformly.
Such a setting is not optimal, considering the curvature of the path. The number of tracing steps was set by default $\# steps = 200$,
and one can ask questions like what is the minimum
number of steps needed or how to set the step-size adaptively (e.g., by setting large tracing steps in regions where the curve is straight,
while finer step-size should be set in neighborhoods where the path is highly curved).
That kind of analysis goes beyond the scope of the current paper.
\end{rem}

%

\section{Numerical examples of derived Gaussian rules}\label{sec:ex}

\subsection{Optimal quadrature rules for $C^0$ uniform quartics, $d=4$, $c=0$}\label{ssec:40}

\begin{figure}[!tb]
\vrule width0pt\hfill
    \begin{overpic}[width=.89\textwidth,angle=0]{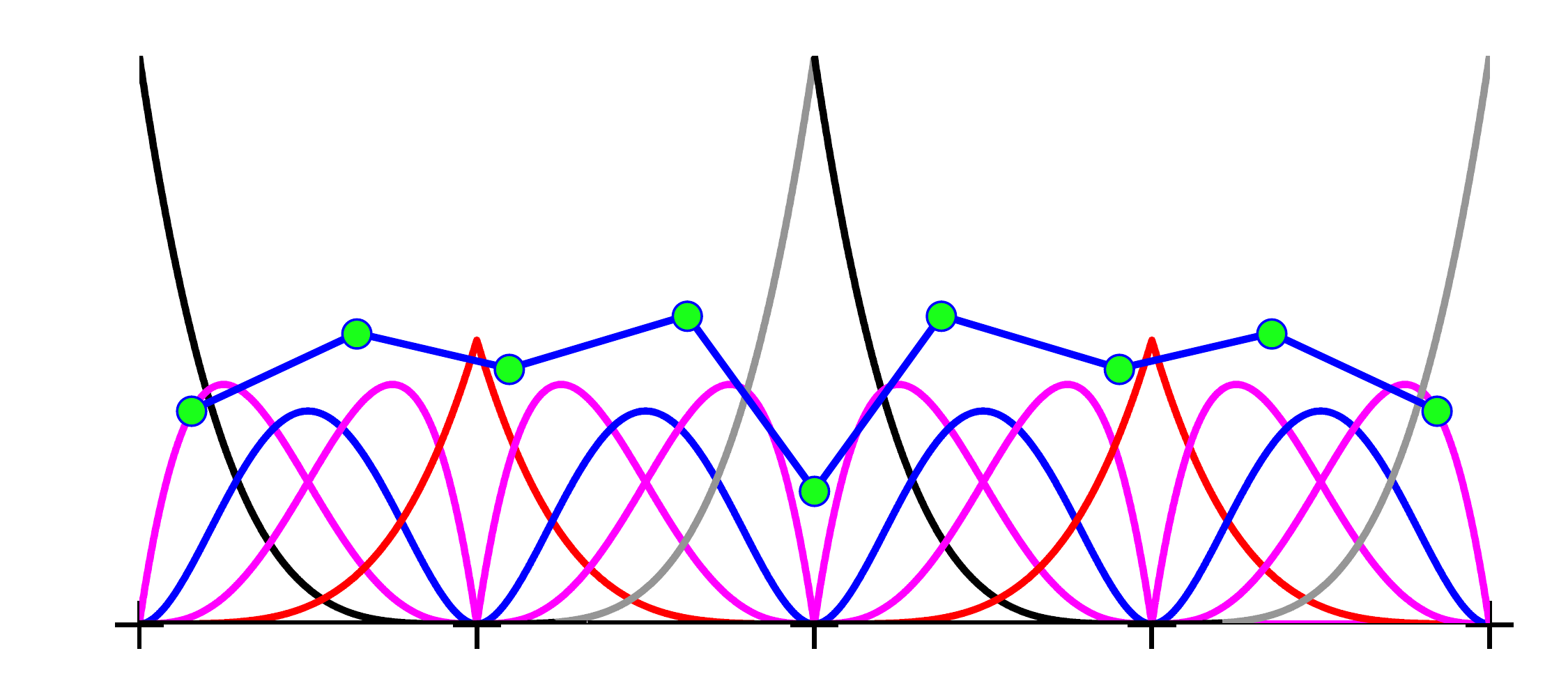}
        \put(11,35){\small$\Dt_{1}$}
        \put(46,35){\small$\Dt_{9}$}
        \put(26,24){\small$\Dt_{5}$}
        \put(25,40){\fcolorbox{gray}{white}{\small$d=4$}}
        \put(40,27){\small$[\tau_4,\omega_4]$}
        \put(2,20){\small$[\tau_1,\omega_1]$}
        \put(5,6){\small$\xt_{0}$}
        \put(96,6){\small$\xt_{4}$}
        \put(7,0){\small $5$}
        \put(28,0){\small $4$}
        \put(52,0){\small $5$}
        \put(73,0){\small $4$}
        \put(95,0){\small $5$}
    \end{overpic}
\hfill \vrule width0pt\\[-4ex]
\Acaption{1ex}{Gauss-Radau source rule. An initial building block for a quartic spline space over $n=4$ normalized elements.
The spline basis functions, and therefore the quadrature rule (green dots), are symmetric with respect to the middle of the interval; the integers bellow each
knot represent the knots' multiplicities. The initial rule is obtained by computing a Gauss-Radau rule for the half of the spline space, spanned by $\Dt_1,\dots,\Dt_9$.
This is achieved by solving two coupled algebraic systems \eqref{eq:IniSystem40a} and \eqref{eq:IniSystem40b}.
The nine unknowns are the four nodes and weights, $\tau_1,\dots, \omega_4$, and the weight $\omega_5$.
 }\label{fig:SourceSpace40}
 \end{figure}

In this case, the original spline space is quadratic and $C^1$-continuous, i.e., $p=2$ and $k=1$ in Section~\ref{sec:matrices}.
For $d=4$, $c=0$, the spline space with an arbitrary number of elements is of odd dimension and,
therefore, we aim at quadrature rules of Gauss-Radau type, requiring one node to be a boundary point or a pre-selected interior point. 

%
%

We initialize the homotopic setting by two Gauss-Radau initial blocks, each of them acting over two elements, see Fig.~\ref{fig:SourceSpace40}.
The spline space is built above a normalized uniform knot vector $\xxt=(0,1,2,3,4)$ with the vector of multiplicities $\bmt=(5,4,5,4,5)$,
that is, the continuities between the elements are in turn $C^0$, $C^{-1}$ and $C^0$. The dimension of the whole space is $18$, but
due to symmetry, we can consider only half of the space because of the discontinuity between the blocks.

The rule is computed by solving two coupled algebraic systems
\begin{equation}\label{eq:IniSystem40a}
\renewcommand{\arraystretch}{1.4}
\begin{array}{rcrcc}
             \omega_1 (1-\tau_1)^4    & + &            \omega_2 (1-\tau_2)^4     & = & \frac{1}{5},  \\
  4 \tau_1   \omega_1 (1-\tau_1)^3    & + & 4 \tau_2   \omega_2 (1-\tau_2)^3     & = & \frac{1}{5},  \\
  6 \tau_1^2 \omega_1 (1-\tau_1)^2    & + & 6 \tau_2^2 \omega_2 (1-\tau_2)^2     & = & \frac{1}{5},  \\
  4 \tau_1^3 \omega_1 (1-\tau_1)^{\sg}& + & 4 \tau_2^3 \omega_2 (1-\tau_2)^{\sg} & = & \frac{1}{5},
\end{array}
\end{equation}
and
\begin{equation}\label{eq:IniSystem40b}
\arraycolsep=1pt
\renewcommand{\arraystretch}{1.4}
\begin{array}{rclcrclrl}
                     &  \omega_3 & (2-\tau_3)^4    & + &                     & \omega_4 & (2-\tau_4)^4    & = & \frac{1}{5} - \rho,  \\
  4 (\tau_3-1)^{\sg} &  \omega_3 & (2-\tau_3)^3    & + & 4 (\tau_4-1)^{\sg}  & \omega_4 & (2-\tau_4)^3    & = & \frac{1}{5},  \\
  6 (\tau_3-1)^2     &  \omega_3 & (2-\tau_3)^2    & + & 6 (\tau_4-1)^2      & \omega_4 & (2-\tau_4)^2    & = & \frac{1}{5},  \\
  4 (\tau_3-1)^3     &  \omega_3 & (2-\tau_3)^{\sg}& + & 4 (\tau_4-1)^3      & \omega_4 & (2-\tau_4)^{\sg}& = & \frac{1}{5},\\
    (\tau_3-1)^4     &  \omega_3 &                 & + &   (\tau_4-1)^4      & \omega_4 &       &+ \sg \frac{1}{2}\omega_5 = & \frac{1}{5},
\end{array}
\end{equation}
where $\rho$ is the residuum of the quadrature rule when applied on $\Dt_5$ on the first element, i.e.,
\begin{equation}\label{eq:residuum}
\rho = \omega_1 \tau_1^4 + \omega_2 \tau_2^4.
\end{equation}
One needs to satisfy nine exactness constraints on the rule when applied to the basis functions. However, we split a large ($9\times 9$)
system into two smaller ones that are coupled via the residuum term of $\Dt_5$,
the only basis function that has support on both elements, see Fig.~\ref{fig:SourceSpace40}.
Factorizing \eqref{eq:IniSystem40a} using computer algebra yields a univariate quadratic polynomial in $\tau_1$ and we obtain
\begin{equation}\label{eq:tau12}
\renewcommand{\arraystretch}{1.4}
\begin{array}{cclcccl}
\tau_{1,2}  & = & \frac{2}{5} \mp \frac{\sqrt{6}}{10},  & \quad & w_{1,2} & = &\frac{4}{9} \mp \frac{\sqrt{6}}{36}.
\end{array}
\end{equation}
Inserting \eqref{eq:tau12} into \eqref{eq:residuum} we obtain the residuum
\begin{equation}\label{eq:residuum2}
\rho = \frac{4}{45}
\end{equation}
and sequentially solving \eqref{eq:IniSystem40b}, we get
\begin{equation*}
\renewcommand{\arraystretch}{1.4}
\begin{array}{cclccclcc}
\tau_{3,4}  & = & \frac{34}{25} \mp \frac{\sqrt{174}}{50},  & \quad & w_{3,4} & = &\frac{76}{153} \mp \frac{7\sqrt{174}}{1972}  & \quad & w_{5} = \frac{4}{17},
\end{array}
\end{equation*}
which, due to symmetry, completes the optimal rule for the considered space.

\begin{figure}[!tb]
 \vrule width0pt\hfill
 \begin{overpic}[width=.89\columnwidth,angle=0]{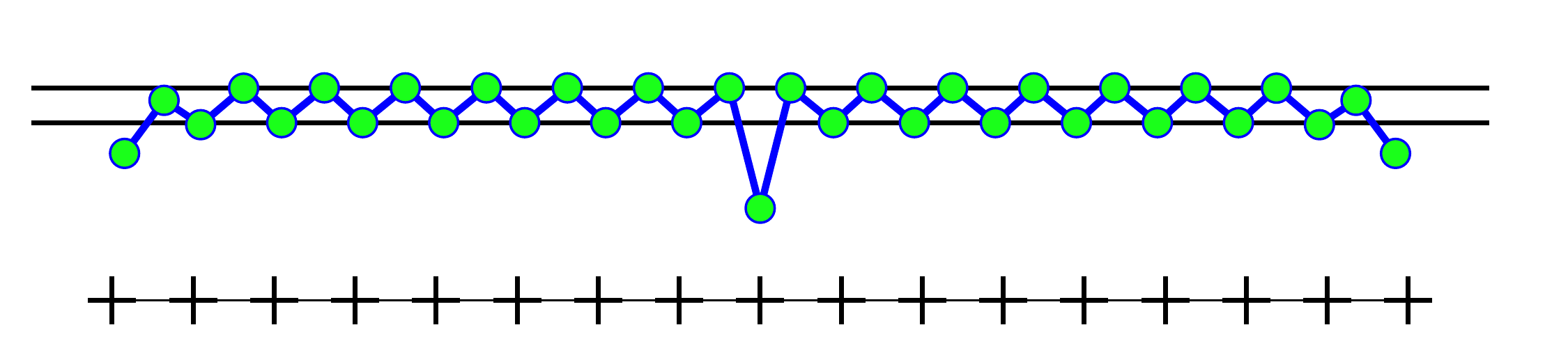}
    \put(11,8){\fcolorbox{gray}{white}{$d=4$, $c=0$, $N=16$}}
    \put(5,0){\footnotesize $0$}
    \put(43,0){\footnotesize $x_M=8$}
    \put(90,0){\footnotesize $x_{16}=16$}
    \put(96,17){\footnotesize $\omega_1^A$}
    \put(96,12){\footnotesize $\omega_2^A$}
	\end{overpic}
 \hfill \vrule width0pt\\[0ex]
 \vrule width0pt\hfill
 \begin{overpic}[width=.99\columnwidth,angle=0]{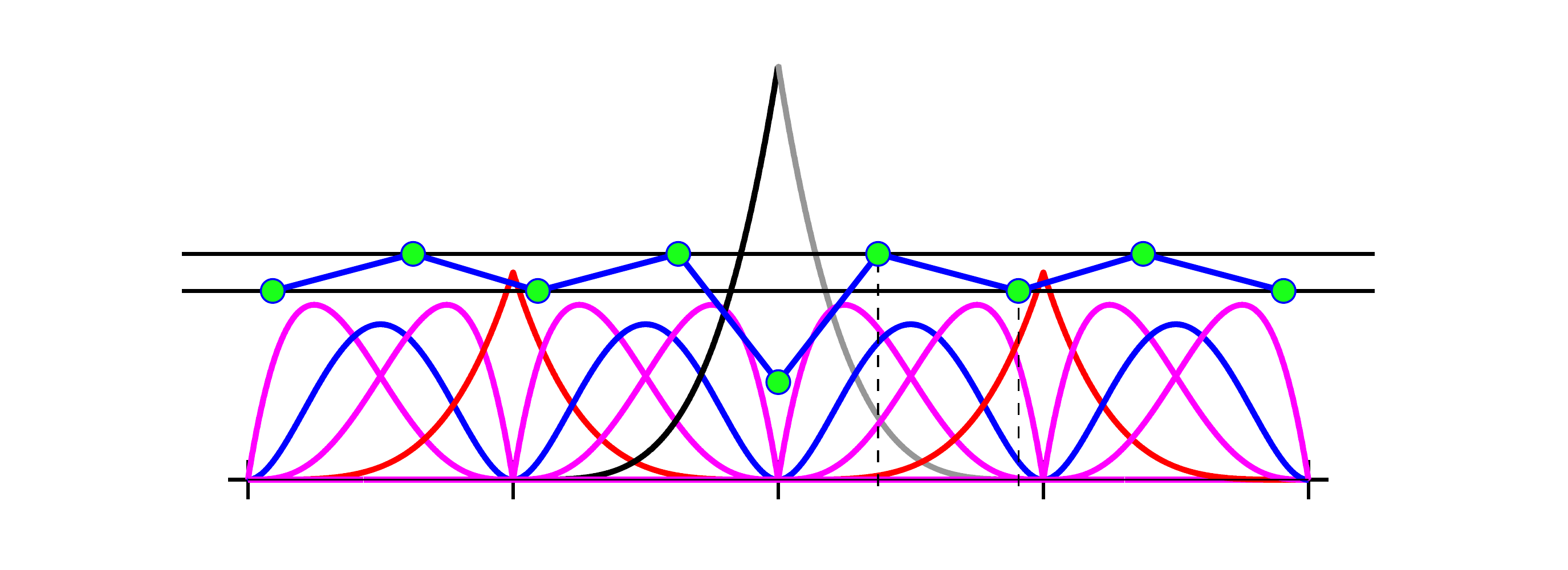}
 \put(15,30){\fcolorbox{gray}{white}{$d=4$, $c=0$, $N=\infty$}}
    \put(60,0){$d_1$}
    \put(49,2){$x_{M}$}
    \put(52,25){$D_{M}$}
    \put(67,2){$x_{M+1}$}
    \put(56,5){$\scriptsize\underbrace{\makebox[1.3cm]{}}$}
    \put(88,21){\footnotesize $\omega_1^A\doteq 0.544$}
    \put(88,16){\footnotesize $\omega_2^A\doteq 0.455$}
	\end{overpic}
 \hfill \vrule width0pt\\[-4ex]
 \Acaption{1ex}{Two-node-per-element rule ($d=4$, $c=0$).
    Top: the layout of the optimal quadrature rule (green dots)
    for uniform knot distribution with $N=16$ elements over $[a,b]=[0,16]$ is shown. 
    All internal knots have multiplicity four, except for the middle one, $x_M$, with multiplicity five.
    Bottom: Asymptotic layout of the optimal rule.
    Two types of basis functions have support only on one element (magenta and blue), while another type spans two (red).
    On the middle pair of elements, the two basis functions that are discontinuous at $x_M$ are shown in black.
    The positions of the nodes are determined by $d_1$ and $d_2$ and, together with the asymptotic weights $\omega_1^A$ and $\omega_2^A$,
    are computed from the $(4\times4)$ asymptotic system \eqref{eq:AsympSystem40}.
    The middle weight is computed from \eqref{eq:OmegaMiddle}, which corresponds to the exactness of the rule when applied on the middle discontinuous
    basis function $D_M$.
 \vspace{0.5cm}
 }\label{fig:40}
 \end{figure}

With this four-element optimal quadrature rule as an elementary quadrature block, we enter the homotopy continuation recursion, see Section \ref{ssec:source},
and derive optimal rules for spaces with various numbers of elements.
In the simplest setup when the target rule is used as a half of the source rule in the next iteration, 
one generates Gauss-Radau rules over $N = 4 \cdot 2^i$, elements in the $i$-th iteration. 
The results of the second iteration of this homotopic setup are shown in Fig.~\ref{fig:40} top, 
where the Gauss-Radau type rule over $16$ uniform elements is derived.

However, one can combine the rules from different levels of recursion, as well as the classical 
Gauss-Radau polynomial rule, to design the source space with the desired number of elements.  
We experimented with this approach and derived the rules with similar errors like the rule showed in Table~\ref{tab:40}. 

\begin{rem}\label{rem:GR}
Note that the rule in Fig.~\ref{fig:40} top is optimal for the space over $N=16$ elements
where all the internal knots have multiplicity four, but the middle knot has multiplicity five.
Consequently, such a rule is also exact for the contained space with all interior knots of multiplicity four (our desired space with $C^0$-continuity
at all knots).
Due to the multiplicity (five) of the middle knot, however, 
one can also interpret the derived rule as two Gauss-Radau rules over $N=8$ elements merged in $C^{-1}$ fashion. 
%
\end{rem}

As $N\rightarrow \infty$, the rule converges to its asymptotic counterpart, an analogy of the midpoint rule of Hughes et al.\cite{Hughes-2010}.
The asymptotic rule for this space is computed from
 \begin{equation}\label{eq:AsympSystem40}
\renewcommand{\arraystretch}{1.4}
\begin{array}{rcrcc}
  4 d_1   \omega_1^A (1-d_1)^3    & + & 4 d_2   \omega_2^A (1-d_2)^3     & = & \frac{1}{5},  \\
  6 d_1^2 \omega_1^A (1-d_1)^2    & + & 6 d_2^2 \omega_2^A (1-d_2)^2     & = & \frac{1}{5},  \\
  4 d_1^3 \omega_1^A (1-d_1)^{\sg}& + & 4 d_2^3 \omega_2^A (1-d_2)^{\sg} & = & \frac{1}{5},  \\
 \omega_1^A d_1^4 + \omega_2^A d_2^4 + d_1 \omega_1^A (1-d_1)^4    & + & d_2   \omega_2^A (1-d_2)^4 & = & \frac{2}{5},
\end{array}
\end{equation}
with four unknowns $d_1$, $d_2$, $\omega_1^A$ and $\omega_2^A$.
This system expresses the exactness of the rule when applied to four consecutive basis functions on a normalized interval
$[x_M,x_{M+1}]=[0,1]$, see Fig.~\ref{fig:40} bottom. The asymptotic middle weight $\omega_M^A$ is computed from
 \begin{equation}\label{eq:OmegaMiddle}
\omega_1^A d_1^4 + \omega_2^A d_2^4 + \frac{1}{2}\omega_M^A = \frac{1}{5}
\end{equation}
which satisfies the exactness of the rule when applied to the middle discontinuous basis function $D_M$.
Solving \eqref{eq:AsympSystem40} with computer algebra software, we obtain the asymptotic values
\begin{equation}\label{eq:AsympVals40}
\renewcommand{\arraystretch}{1.4}
\begin{array}{cclcc}
d_{1}       & = & \frac{1}{2} + \frac{\sqrt{7}}{10} - \frac{\sqrt{2}}{10} & \doteq & 0.62315377486914955417, \\
d_{2}       & = & \frac{1}{2} - \frac{\sqrt{7}}{10} - \frac{\sqrt{2}}{10} & \doteq & 0.09400351265623143607, \\
\omega_1^A  & = & \frac{1}{2} + \frac{\sqrt{14}}{84}                      & \doteq & 0.54454354031873739745, \\
\omega_2^A  & = & \frac{1}{2} - \frac{\sqrt{14}}{84}                      & \doteq & 0.45545645968126260255, \\
\omega_M^A  & = & \frac{\sqrt{2}}{6}                                      & \doteq & 0.23570226039551584147 \\
\end{array}
\end{equation}
which define the optimal asymptotic rule 
\begin{equation}\label{eq:limit40}
\begin{split}
\int_{\mathbb R}f(t) = \omega_M^A f(x_M)  & + \sum\limits_{i \in \mathbb Z^+} h(\omega_1^A f((i-d_1)h) + \omega_2^A f((i-d_2)h)) \\
                                          & + \sum\limits_{i \in \mathbb Z^-} h(\omega_1^A f((i+d_1)h) + \omega_2^A f((i+d_2)h)).
\end{split}
\end{equation}

\begin{table}[!tb]
 \begin{center}
  \begin{minipage}{0.9\textwidth}
\caption{Two-node patterned Gaussian quadrature rule
for $d=4$, $c=0$, with $N=32$ uniform elements over $[0,N]$. Observe the convergence to the asymptotic values \eqref{eq:AsympVals40}.
The nodes and weights are shown with the precision of $20$ decimal digits and the values on the first ten boundary elements differ
from the asymptotic values by more than $16$ decimal digits.
}\label{tab:40}
  \end{minipage}
\vspace{0.2cm}\\
\small{
\renewcommand{\arraystretch}{1.15}
\renewcommand{\tf}{\small}
\begin{tabular}{| c | c || r| r|}\hline
\multicolumn{2}{|c}{} & \multicolumn{2}{c|}{$d=4$, $c=0$, $N=32$, uniform, $\|\br\| = 4.81^{-26}$ } \\\hline
$\#$el. & $i$ &  \multicolumn{1}{c|}{$\tau_i$}  & \multicolumn{1}{c|}{$\omega_i$} \\\hline\hline
\multirow{2}{*}{1} & 1 & \tf 0.15505102572168219018 & \tf 0.37640306270046727505   \\
                   & 2 & \tf 0.64494897427831780982 & \tf 0.51248582618842161384   \\\hline
\multirow{2}{*}{2} & 3 & \tf 1.09618188083454161658 & \tf 0.44990832345215269846   \\
                   & 4 & \tf 1.62381811916545838342 & \tf 0.54355572883542900089   \\\hline
\multirow{2}{*}{3} & 5 & \tf 2.09406803063701196217 & \tf 0.45528750742625502979   \\
                   & 6 & \tf 2.62317334867333286542 & \tf 0.54451443355215653685   \\\hline
\multirow{2}{*}{4} & 7 & \tf 3.09400541223051380344 & \tf 0.45545148116758058646   \\
                   & 8 & \tf 3.62315435108309566402 & \tf 0.54454268347129809054   \\\hline
\multirow{2}{*}{5} & 9 & \tf 4.09400356857477400144 & \tf 0.45545631312314607882   \\
                   & 10& \tf 4.62315379183131736912 & \tf 0.54454351509548282068   \\\hline
\multirow{2}{*}{6} & 11& \tf 5.09400351430231989540 & \tf 0.45545645536699068802   \\
                   & 12& \tf 5.62315377536846916574 & \tf 0.54454353957623409006   \\\hline
\multirow{2}{*}{7} & 13& \tf 6.09400351270468775630 & \tf 0.45545645955426229179   \\
                   & 14& \tf 6.62315377488384815118 & \tf 0.54454354029688014082   \\\hline
\multirow{2}{*}{8} & 15& \tf 7.09400351265765785696 & \tf 0.45545645967752406246   \\
                   & 16& \tf 7.62315377486958224046 & \tf 0.54454354031809397989   \\\hline
\multirow{2}{*}{9} & 17& \tf 8.09400351265627342598 & \tf 0.45545645968115255021   \\
                   & 18& \tf 8.62315377486916229127 & \tf 0.54454354031871845700   \\\hline
\multirow{2}{*}{10}& 19& \tf 9.09400351265623267214 & \tf 0.45545645968125936291   \\
                   & 20& \tf 9.62315377486914992912 & \tf 0.54454354031873683989   \\\hline
\multirow{2}{*}{11}& 21& \tf10.09400351265623147246 & \tf 0.45545645968126250719   \\
                   & 22& \tf10.62315377486914956521 & \tf 0.54454354031873738103   \\\hline
\multirow{2}{*}{12}& 23& \tf11.09400351265623143714 & \tf 0.45545645968126259975   \\
                   & 24& \tf11.62315377486914955449 & \tf 0.54454354031873739696   \\\hline
\multirow{2}{*}{13}& 25& \tf12.09400351265623143610 & \tf 0.45545645968126260247   \\
                   & 26& \tf12.62315377486914955418 & \tf 0.54454354031873739743   \\\hline
\multirow{2}{*}{14}& 27& \tf13.09400351265623143607 & \tf 0.45545645968126260255   \\
                   & 28& \tf13.62315377486914955417 & \tf 0.54454354031873739745   \\\hline
\multirow{2}{*}{15}& 29& \tf14.09400351265623143607 & \tf 0.45545645968126260255   \\
                   & 30& \tf14.62315377486914955417 & \tf 0.54454354031873739745   \\\hline
\multirow{3}{*}{16}& 31& \tf15.09400351265623143607 & \tf 0.45545645968126260255   \\
                   & 32& \tf15.62315377486914955417 & \tf 0.54454354031873739745   \\
                   & 33& \multicolumn{1}{l|}{$\tf16$}&\tf 0.23570226039551584147   \\\hline
\end{tabular}
}
\end{center}
\end{table}

Table~\ref{tab:40} shows the rule derived by our algorithm for $N=32$ and we observe convergence to the asymptotic rule.
The first ten nodes and weights differ from the asymptotic values. The optimal rule at hand
for $N>32$ consists of the ten-elements quadrature block from Table~\ref{tab:40} (lines $1$ to $20$)
and the analytic asymptotic rule with nodes and weights \eqref{eq:AsympVals40}.

The error of the rule $\Q$ is expressed in terms of the Euclidean norm of the vector of the residues of the
system (\ref{eq:IniSystem}), normalized by the dimension of the system
\begin{equation}\label{eq:Error}
\|\br\| = \frac{1}{2m}(\sum_{i=1}^{2m} (\Q_a^b[D_i] - I[D_i])^2)^{\frac{1}{2}}.
\end{equation}

\subsection{Optimal rules for $C^1$ uniform sixtics, $d=6$, $c=1$}\label{ssec:61}

For this space, the dimension is even for an arbitrary even number of elements and therefore admits
an optimal rule without forcing additional constraints, in contrast to Section \ref{ssec:40}.
We initialize the homotopic setup by the optimal rule for a two-element building block ($N=2$), see \eqref{eq:NW61}
and Fig.~\ref{fig:BasicBlock61}. The homotopic process deriving the optimal rule from two building blocks ($N=4$)
is shown in Fig.~\ref{fig:Trans61}. Further, we proceed recursively as described in Section~\ref{ssec:source}.

 \begin{figure}[!tb]
 \vrule width0pt\hfill
 \begin{overpic}[width=.89\columnwidth,angle=0]{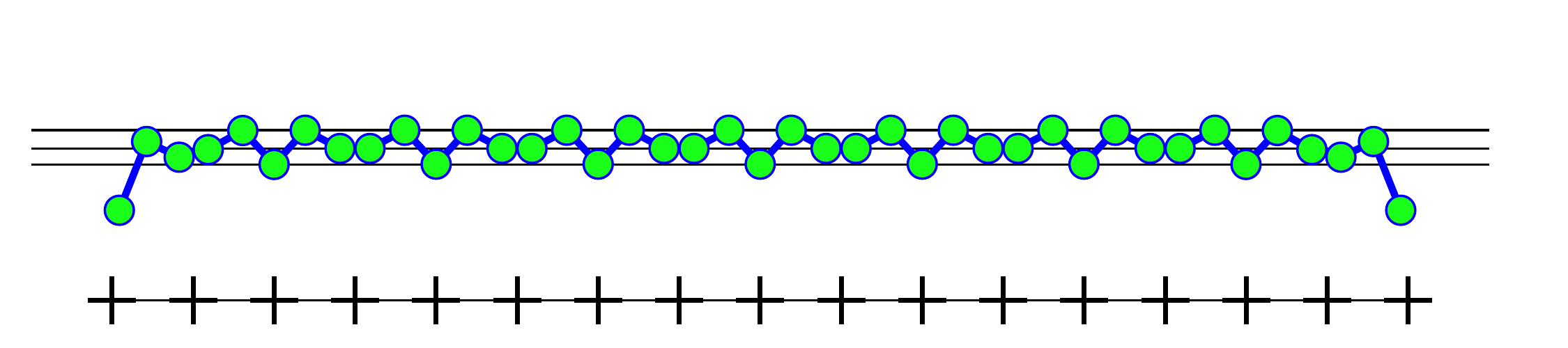}
    \put(11,7){\fcolorbox{gray}{white}{$d=6$, $c=1$, $N=16$}}
    \put(5,0){\footnotesize $0$}
    \put(90,0){\footnotesize $x_{16}=16$}
	\end{overpic}
 \hfill \vrule width0pt\\[2ex]
 \vrule width0pt\hfill
 \begin{overpic}[width=.99\columnwidth,angle=0]{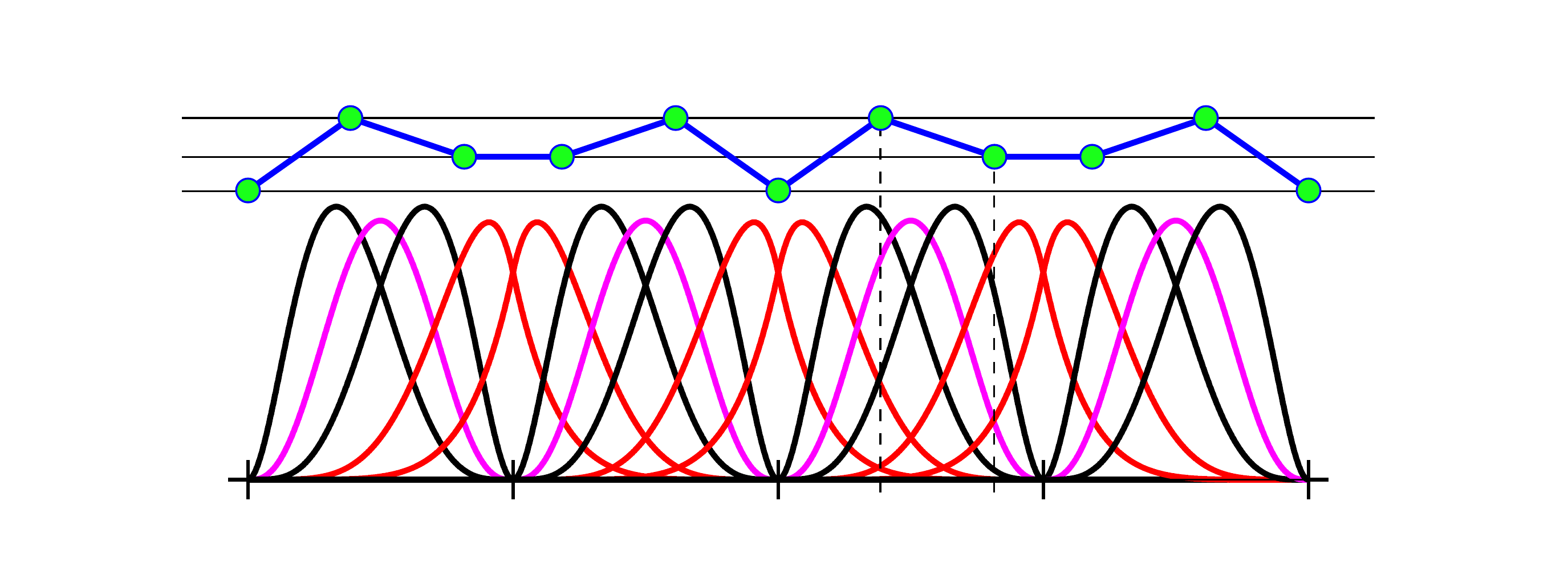}
 \put(15,33){\fcolorbox{gray}{white}{$d=6$, $c=1$, $N=\infty$}}
    \put(60,0){$d_1$}
    \put(15,2){$x_{i-2}$}
    \put(49,2){$x_{i}$}
    \put(67,2){$x_{i+1}$}
    \put(56,5){$\scriptsize\underbrace{\makebox[1.3cm]{}}$}
    \put(88,30){\footnotesize $\omega_1^A\doteq 0.436$}
    \put(88,26){\footnotesize $\omega_2^A\doteq 0.389$}
    \put(88,22){\footnotesize $\omega_3^A\doteq 0.349$}
	\end{overpic}
 \hfill \vrule width0pt\\[-4ex]
 \Acaption{1ex}{Two-and-half point rule ($d=6$, $c=1$).
    Top: the optimal quadrature rule (green dots)
    for uniform knot distribution with $N=16$ elements over $[a,b]=[0,16]$ is shown.
    For growing $N$, the rule quickly converges to the asymptotic rule (bottom)
    where one set of nodes become even knots. The other two sets are determined by $d_1$ and
    $d_2$ and together with the asymptotic weights $\omega_j^A$, $j=1,2,3$  are computed from
    the system \eqref{eq:AsympSystem61}.
 \vspace{0.5cm}
 }\label{fig:61}
 \end{figure}

The rule for $N=16$ is shown in Fig.~\ref{fig:61} and we again observe fast convergence to the asymptotic counterpart ($N=\infty$).
The asymptotic rule requires five nodes every two elements and the layout of the nodes is shown in Fig.~\ref{fig:61} bottom.
The asymptotic rule is computed from
\begin{equation}\label{eq:AsympSystem61}
\renewcommand{\arraystretch}{1.4}
\begin{array}{rcrcc}
 15 d_1^2 \omega_1^A (1-d_1)^4 & + &15 d_2^2 \omega_2^A (1-d_2)^4 & = & \frac{1}{7},  \\
 20 d_1^3 \omega_1^A (1-d_1)^3 & + &20 d_2^3 \omega_2^A (1-d_2)^3 & = & \frac{1}{7},  \\
 15 d_1^4 \omega_1^A (1-d_1)^2 & + &15 d_2^4 \omega_2^A (1-d_2)^2 & = & \frac{1}{7},  \\
 6 \omega_1^A d_1^5 -5\omega_1^A d_1^6 & + & 6 \omega_2^A d_2^5 -5\omega_2^A d_2^6 & = & \frac{2}{7},
\end{array}
\end{equation}
that expresses exactness of the rule when applied to four consecutive basis functions,
three having support on only one element whilst the fourth one spans two.
%
The four unknowns are $d_1$ and $d_2$ (that determine the positions of the nodes) and the asymptotic weights $\omega_1^A$ and $\omega_2^A$.
The last asymptotic weight, $\omega_3^A$, is sequentially computed from
\begin{equation}\label{eq:W3}
2\omega_1^A + 2\omega_2^A + \omega_3^A = 2
\end{equation}
which is the exactness of the rule when applied to a constant function.
The algebraic factorization of \eqref{eq:AsympSystem61} gives
\begin{equation}\label{eq:Uni61Asymp}
52-364d_1 + 905d_1^2-938d_1^3+343d_1^4,
\end{equation}
which, due to a quartic degree only, admits a closed-form formula for the asymptotic nodes and weights
%
\begin{equation}\label{eq:AsympVals61}
\renewcommand{\arraystretch}{1.5}
\begin{array}{cclcc}
d_{1}       & = & \frac{67}{98} - \frac{3\sqrt{78}}{98} - \frac{\sqrt{95-10\sqrt{78}}}{98} & \doteq & 0.38693556354866909100, \\
d_{2}       & = & \frac{67}{98} + \frac{3\sqrt{78}}{98} - \frac{\sqrt{95+10\sqrt{78}}}{98} & \doteq & 0.81587550281258499773, \\
\omega_1^A  & = & \frac{1693}{4160} + \frac{3\sqrt{5}\sqrt{13}}{4160}
                - \frac{673\sqrt{5}\sqrt{6}}{99840} + \frac{2047\sqrt{6}\sqrt{13}}{299520} & \doteq & 0.43622310273429582467, \\
\omega_2^A  & = & \frac{1693}{4160} + \frac{3\sqrt{5}\sqrt{13}}{4160}
                + \frac{673\sqrt{5}\sqrt{6}}{99840} - \frac{2047\sqrt{6}\sqrt{13}}{299520} & \doteq & 0.38934746132575016040, \\
\omega_3^A  & = & \frac{387}{1040} - \frac{3\sqrt{5}\sqrt{13}}{1040}                       & \doteq & 0.34885887187990802985. \\
\end{array}
\end{equation}
%
Finally, the asymptotic rule for an infinite uniform knot vector with the elements of size $h$ is formalized as
\begin{equation}\label{eq:limit61}
\begin{split}
\int_{\mathbb R}f(t) = \sum\limits_{i \in \mathbb Z} h( &  \omega_1^A (f((2i+d_1)h) + f((2i+2-d_1)h))  \\
                                                      + &  \omega_2^A (f((2i+d_2)h) + f((2i+2-d_2)h)) + \omega_3^A f(2ih) )
\end{split}
\end{equation}
which is exact and optimal over the real line. The convergence of the rules over finite domains derived
via homotopy continuation is shown in Table~\ref{tab:61}. For $N=16$, only the nodes and weights on the first five boundary elements
differ from the asymptotic values \eqref{eq:AsympVals61}. Thus the optimal rule for finite domains
is therefore a combination of the nodes and weights on the first five boundary elements (lines $1$ to $13$ in Table~\ref{tab:61})
and the asymptotic rule \eqref{eq:limit61} for all the intermediate elements.

We emphasize here that our optimal rules require (asymptotically when $N$ is large) only $2.5$ nodes per element 
in contrast to $4$ nodes required by classical polynomial Gauss rule. This reduction might not be that significant in one variable, however, 
for 3D problems when tensor product rules are used the reduction ratio is already $(\frac{2.5}{4})^3 \doteq 24\%$.

\begin{table}[!tb]
\begin{center}
  \begin{minipage}{0.9\textwidth}
\caption{Gaussian quadrature rule the Galerkin space $d=6$, $c=1$, with $N=16$ uniform elements over $[0,N]$.
Observe fast convergence to the asymptotic rule, see \eqref{eq:AsympVals61};
the nodes and weights only on the first five boundary elements differ from the asymptotic values by more than $16$ decimal digits.
}\label{tab:61}
  \end{minipage}
\vspace{0.2cm}\\
\small{
\renewcommand{\arraystretch}{1.15}
\renewcommand{\tf}{\small}
\begin{tabular}{| c | c || l| l|}\hline
\multicolumn{2}{|c}{} & \multicolumn{2}{c|}{$d=6$, $c=1$, $N=16$, uniform, $\|\br\| = 3.75^{-26}$ } \\\hline
$\#$el. & $i$ &  \multicolumn{1}{c|}{$\tau_i$}  & \multicolumn{1}{c|}{$\omega_i$} \\\hline\hline
\multirow{3}{*}{1} & 1 & \tf 0.09260767873646902812 & \tf 0.23050486991521396993   \\
                   & 2 & \tf 0.42847197760814208611 & \tf 0.40704416177654188371   \\
                   & 3 & \tf 0.83018935543014295850 & \tf 0.36711516474717107854   \\\hline
\multirow{2}{*}{2} & 4 & \tf 1.18644180845680657718 & \tf 0.38605131464693100757   \\
                   & 5 & \tf 1.61390002454892326539 & \tf 0.43521953213902864887   \\\hline
\multirow{3}{*}{3} & 6 & \tf 2.00010871499078850047 & \tf 0.34849458018527149253   \\
                   & 7 & \tf 2.38693570464281488360 & \tf 0.43622300768518266759   \\
                   & 8 & \tf 2.81587555220352588540 & \tf 0.38934738499907207358   \\\hline
\multirow{2}{*}{4} & 9 & \tf 3.18412450505465915622 & \tf 0.38934744984465969166   \\
                   & 10& \tf 3.61306443926733132981 & \tf 0.43622309934864369784   \\\hline
\multirow{3}{*}{5} & 11& \tf 4.00000000036580449734 & \tf 0.34885887065223780524   \\
                   & 12& \tf 4.38693556354866909260 & \tf 0.43622310273429582360   \\
                   & 13& \tf 4.81587550281258499829 & \tf 0.38934746132575015954   \\\hline
\multirow{2}{*}{6} & 14& \tf 5.18412449718741500236 & \tf 0.38934746132575016027   \\
                   & 15& \tf 5.61306443645133090903 & \tf 0.43622310273429582463   \\\hline
\multirow{3}{*}{7} & 16& \tf 6                      & \tf 0.34885887187990802983   \\
                   & 17& \tf 6.38693556354866909100 & \tf 0.43622310273429582467   \\
                   & 18& \tf 6.81587550281258499773 & \tf 0.38934746132575016040   \\\hline
\multirow{2}{*}{8} & 19& \tf 7.18412449718741500227 & \tf 0.38934746132575016040   \\
                   & 20& \tf 7.61306443645133090900 & \tf 0.43622310273429582467   \\
                   & 21& \tf 8                      & \tf 0.34885887187990802984   \\\hline
\end{tabular}
}
\end{center}
\end{table}

\begin{figure}[!tb]
 \vrule width0pt\hfill
 \begin{overpic}[width=.89\columnwidth,angle=0]{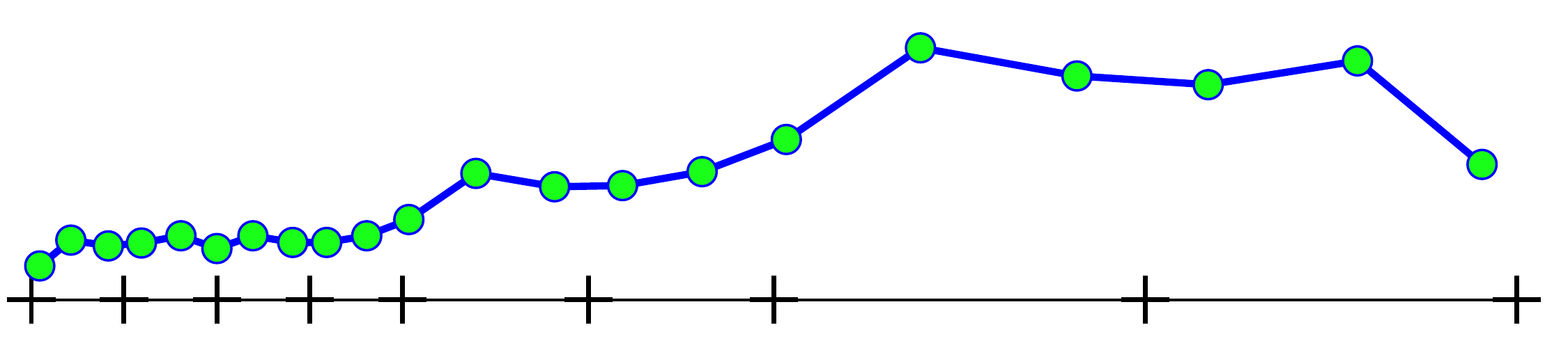}
    \put(-2,18){\fcolorbox{gray}{white}{$d=6$, $c=1$, $N=8$, non-uniform}}
    \put(0,0){\footnotesize $0$}
    \put(93,0){\footnotesize $x_{8}=8$}
	\end{overpic}
 \hfill \vrule width0pt\\[-4ex]
 \Acaption{1ex}{Gaussian rule (green dots) for non-uniform sixtic spline space with $N=8$
 elements over $[a,b]=[0,8]$ is shown.
    The partition knot vector is $\bx=(0,\frac{1}{2},1,\frac{3}{2},2,3,4,6,8)$
    where all the internal knot have multiplicity five ($c=1$) and the boundary knots
    have multiplicity seven (open end condition). The nodes and weights are listed in
    Table \ref{tab:61NonUni}.
    }\label{fig:61NonUni}
 \end{figure}

While the uniform knot sequences converge to the asymptotic counterparts and therefore only several boundary nodes
and weights differ from the regular pattern, this simple behavior is no longer available for general knot sequences. The homotopic concept, however,
is well suited especially for non-uniform spline spaces. An example of a Gaussian rule for such a space is shown in Fig.~\ref{fig:61NonUni} and
the nodes and weights are listed in Table~\ref{tab:61NonUni}. The $C^1$-continuous sixtic spline space spans eight elements and the knot distribution is denser closer
to the left boundary. 
For applications like boundary layer flows \cite{Bazilevs-2007-Variational,Motlagh-2013-Simulation,Ayotte-2004-BoundLayer} where a finer mesh resolution is typically needed close to one boundary, our homotopic approach offers a tool to derive optimal rules for these non-uniform spaces.

\section{Conclusion}\label{sec:conl}

We derive optimal (Gaussian) quadrature rules for spline spaces frequently appearing
in Galerkin discretizations where the original spaces are $C^1$ quadratics or $C^2$ cubics.
The rules are optimal, that is, they require the minimum possible quadrature points, and
are exact over finite intervals with uniformly distributed elements.
We have numerically shown that the presented rules quickly converge
to their asymptotic counterparts, as the nodes and weights on only few boundary elements
differ from the repetitive asymptotic pattern.
Due to the analytical nature of the asymptotic rules, we believe our rules can be easily incorporated to IGA software libraries
like PetIGA \cite{petsc,petiga}, as only few boundary data must be read from a look-up table.

Moreover, all the derived rules possess positive weights which make them a convenient choice from the
point of view of numerical stability, e.g., when compared to \cite{Hughes-2012}
where some weights are close to zero and are negative.
We have also shown that the homotopic concept is not limited to uniform knot spans
and the rules can be adapted to non-uniform spacing.

\section*{Acknowledgements}

This publication was made possible in part by a National Priorities Research Program grant 7-1482-1-278 from the Qatar National Research Fund (a member of The Qatar Foundation), by the European Union’s Horizon 2020 Research and Innovation Program under the Marie Sklodowska-Curie grant agreement No. 644602, and the Center for Numerical Porous Media at King Abdullah University of Science and Technology (KAUST). The J. Tinsley Oden Faculty Fellowship Research Program at the Institute for Computational Engineering and Sciences (ICES) of the University of Texas at Austin has partially supported the visits of VMC to ICES.

 \begin{table}[!tb]
\begin{center}
  \begin{minipage}{0.9\textwidth}
\caption{Gaussian quadrature rule for a non-uniform spline space, $d=6$, $c=1$, with $N=8$ elements over $[0,N]$
shown in Fig.~\ref{fig:61NonUni}. The dimension of the space is $42$ and therefore the optimal rule requires $21$ nodes.
}\label{tab:61NonUni}
  \end{minipage}
\vspace{0.2cm}\\
\small{
\renewcommand{\arraystretch}{1.15}
\renewcommand{\tf}{\small}
\begin{tabular}{| c | c || l| l|}\hline
\multicolumn{2}{|c}{} & \multicolumn{2}{c|}{$d=6$, $c=1$, $N=8$, non-uniform, $\|\br\| = 8.57^{-30}$ } \\\hline
$\#$el. & $i$ &  \multicolumn{1}{c|}{$\tau_i$}  & \multicolumn{1}{c|}{$\omega_i$} \\\hline\hline
\multirow{3}{*}{1} & 1 & \tf 0.04630383936823451406 & \tf 0.11525243495760698496   \\
                   & 2 & \tf 0.21423598880407104306 & \tf 0.20352208088827094186   \\
                   & 3 & \tf 0.41509467771507147925 & \tf 0.18355758237358553927   \\\hline
\multirow{2}{*}{2} & 4 & \tf 0.59322090422840328859 & \tf 0.19302565732346550379   \\
                   & 5 & \tf 0.80695001227446163269 & \tf 0.21760976606951432444   \\\hline
\multirow{3}{*}{3} & 6 & \tf 1.00005435749539425024 & \tf 0.17424729009263574626   \\
                   & 7 & \tf 1.19346785232140744180 & \tf 0.21811150384259133380   \\
                   & 8 & \tf 1.40793777610176294270 & \tf 0.19467369249953603679   \\\hline
\multirow{2}{*}{4} & 9 & \tf 1.59206225252732957811 & \tf 0.19467372492232984583   \\
                   & 10& \tf 1.80653221963366566491 & \tf 0.21811154967432184892   \\\hline
\multirow{3}{*}{5} & 11& \tf 2.03366386534871873978 & \tf 0.27364402258520424593   \\
                   & 12& \tf 2.39575347568220124424 & \tf 0.42990626936051039389   \\
                   & 13& \tf 2.81890006050280681835 & \tf 0.38464672961950394215   \\\hline
\multirow{2}{*}{6} & 14& \tf 3.18460630101439855425 & \tf 0.38864808057905118797   \\
                   & 15& \tf 3.61323715670019192625 & \tf 0.43601548697564552637   \\\hline
\multirow{3}{*}{7} & 16& \tf 4.06704953147532718337 & \tf 0.54635960217072361337   \\
                   & 17& \tf 4.78975598662033980891 & \tf 0.85789420372567177811   \\
                   & 18& \tf 5.63316509361482355771 & \tf 0.76272937432250973703   \\\hline
\multirow{2}{*}{8} & 19& \tf 6.34055900169025774853 & \tf 0.73283097829499297885   \\
                   & 20& \tf 7.14341666786039006430 & \tf 0.81371802826546978692   \\
                   & 21& \tf 7.81485959249475117486 & \tf 0.46082194145685870291   \\\hline
\end{tabular}
}
\end{center}
\end{table}


\bibliographystyle{plain}
\bibliography{StiffnessMatrixEvenArXiv}

\end{document}